\documentclass[reqno,12pt]{amsart} 

\usepackage{bbm,enumerate,amsmath,amstext,mathrsfs,amssymb,graphicx}

\usepackage[latin1]{inputenc}
\usepackage[T1]{fontenc}
\usepackage[english]{babel}

\newcommand{\R}{\mathbbm{R}}
\newcommand{\N}{\mathbbm{N}}
\newcommand{\Q}{\mathbbm{Q}}
\newcommand{\Z}{\mathbbm{Z}}
\newcommand{\C}{\mathbbm{C}}

\newcommand{\CP}{\C P}

\newcommand{\thh}{\mathrm{Th}}
\newcommand{\hocolim}{\operatorname*{hocolim}}
\newcommand{\Diff}{\mathrm{Diff}}
\newcommand{\Int}{\mathrm{int}}

\newcommand{\Emb}{\mathrm{Emb}}

\newcommand{\st}{\mathrm{st}}

\newcommand{\dist}{\mathrm{dist}}
\newcommand{\Aut}{\mathrm{Aut}}
\newcommand{\Out}{\mathrm{Out}}
\newcommand{\Ver}{\mathrm{Vert}}
\newcommand{\Simp}{\mathrm{Simp}}
\newcommand{\sd}{\mathrm{sd}}
\newcommand{\wreath}{\wr}

\newcommand{\Sub}{\mathrm{Sub}}

\newcommand{\IM}{\mathrm{Im}}

\renewcommand{\phi}{\varphi}
\renewcommand{\epsilon}{\varepsilon}

\newcommand{\DD}{\mathscr{D}}

\newcommand{\CC}{\mathscr{C}}

\newcommand{\id}{\mathrm{id}}

\newcommand{\colim}{\operatorname*{colim}}

\newcommand{\Map}{\mathrm{Map}}
\newcommand{\Gr}{\mathrm{Gr}}

\newcommand{\EE}{\mathscr{E}}
\newcommand{\VV}{\mathscr{V}}
\newcommand{\UU}{\mathscr{U}}

\newcommand{\Ob}{\mathrm{ob}}
\newcommand{\mor}{\mathrm{mor}}

\newcommand{\supp}{\mathrm{supp}}

\newcommand{\MT}[1]{\mathit{MT}O(#1)}

\usepackage{amsthm}
\theoremstyle{plain}
\newtheorem{theorem}{Theorem}[section]
\newtheorem{proposition}[theorem]{Proposition}
\newtheorem{lemma}[theorem]{Lemma}

\newtheorem{corollary}[theorem]{Corollary}

\theoremstyle{definition}
\newtheorem{definition}[theorem]{Definition}

\newtheorem{example}[theorem]{Example}

\theoremstyle{remark}

\newtheorem{construction}[theorem]{Construction}
\newtheorem{remark}[theorem]{Remark}
\newtheorem*{remark*}{Remark}

\numberwithin{equation}{section}

\usepackage[all]{xy}
\CompileMatrices

\begin{document}

\title{Stable Homology of Automorphism groups of Free Groups}
\author{S\o ren Galatius}
\date{\today}
\thanks{S.\ Galatius is supported by NSF grant DMS-0505740}

\begin{abstract}
  Homology of the group $\Aut(F_n)$ of automorphisms of a free group
  on $n$ generators is known to be independent of $n$ in a certain
  stable range.  Using tools from homotopy theory, we prove that in
  this range it agrees with homology of symmetric groups.  In
  particular we confirm the conjecture (\cite{MR1671188}) that stable
  rational homology of $\Aut(F_n)$ vanishes.
\end{abstract}

\maketitle\let\languagename\relax

\tableofcontents

\section{Introduction}
\label{sec:intro}

\subsection{Results}
\label{subsec:results}
Let $F_n = \langle x_1, \dots, x_n\rangle$ be the free group on $n$
generators and let $\Aut(F_n)$ be its automorphism group.  Let
$\Sigma_n$ be the symmetric group and let $\phi_n: \Sigma_n \to
\Aut(F_n)$ be the homomorphism that to a permutation $\sigma$
associates the automorphism $\phi(\sigma): x_i \mapsto x_{\sigma(i)}$.
The main result of the paper is the following theorem.
\begin{theorem}\label{thm:main}
  $\phi_n$ induces an isomorphism
  \begin{align*}
    (\phi_n)_*\colon H_k(\Sigma_n) \to H_k(\Aut(F_n))
  \end{align*}
  for $n > 2k+1$.
\end{theorem}
The homology groups in the theorem are independent of $n$ in the sense
that increasing $n$ induces isomorphisms $H_k(\Sigma_n) \cong
H_k(\Sigma_{n+1})$ and $H_k(\Aut(F_n)) \cong H_k(\Aut(F_{n+1}))$ when
$n > 2k+1$.  For the symmetric group this was proved by Nakaoka
(\cite{MR0112134}) and for $\Aut(F_n)$ by Hatcher and Vogtmann
(\cite{MR1678155},\cite{MR2113904}).  The homology groups
$H_k(\Sigma_n)$ are completely known.  With finite coefficients the
calculation was done by Nakaoka and can be found in \cite{MR0131874}.
We will not quote the result here.  With rational coefficients the
homology groups vanish because $\Sigma_n$ is a finite group, so
theorem~\ref{thm:main} has the following corollary.
\begin{corollary}\label{cor:Mumford}
  The groups
  \begin{align*}
    H_k(\Aut(F_n);\Q)
  \end{align*}
  vanish for $n > 2k+1$.
\end{corollary}

The groups $\Aut(F_n)$ are special cases of a more general series of
groups $A_n^s$, studied in \cite{MR2113904} and
\cite{MR2220689}.  We recall the definition.  For a finite graph
$G$ without vertices of valence $0$ and 2, let $\partial G$ denote the
set of vertices of valence 1.  Let $h\Aut(G)$ denote the topological
monoid of homotopy equivalences $G \to G$ that restrict to the
identity map on $\partial G$.  Let $\Aut(G) = \pi_0 h\Aut(G)$.
\begin{definition}
  Let $G_n^s$ be a connected graph with $s$ leaves and first Betti
  number $b_1(G_n^s) = n$.  For $s + n \geq 2$ let
  \begin{align*}
    A_n^s = \Aut(G_n^s).
  \end{align*}
  In particular $A_n^0 = \Out(F_n)$ and $A_n^1 = \Aut(F_n)$.  $A_0^s$
  is the trivial group for all $s$.
\end{definition}

There are natural group maps for $n \geq 0$, $s \geq 1$
\begin{align}\label{eq:12}
  \xymatrix{
    {A_n^{s-1}} & {A_n^s} \ar[l]_-{\alpha_n^s} \ar[r]^-{\beta_n^s} &
    {A_n^{s+1}} \ar[r]^-{\gamma_n^s} & {A_{n+1}^s}.
  }
\end{align}
$\beta_n^s$ and $\gamma_n^s$ are induced by gluing a $Y$-shaped graph
to $G_n^s$ along part of $\partial G_n^s$. $\alpha_n^s$ is induced by
collapsing a leaf.  We quote the following theorem.
\begin{theorem}[\cite{MR2113904},\cite{MR2220689}]
  \label{thm:HVW}
  $(\beta_n^s)_*$ and $(\gamma_n^s)_*$ are isomorphisms for $n >
  2k+1$.  $(\alpha_n^s)_*$ is an isomorphism for $n > 2k+1$ for $s >
  1$ and $(\alpha_n^1)_*$ is an isomorphism for $n > 2k+3$.
\end{theorem}

The main theorem~\ref{thm:main} calculates the homology of these
groups in the range in which it is independent of $n$ and $s$.  In
other words, we calculate the homology of the group
\begin{align*}
  \Aut_\infty = \colim_{n\to \infty} \Aut(F_n).
\end{align*}
An equivalent formulation of the main theorem is that the map of
classifying spaces $B\Sigma_\infty \to B\Aut_\infty$ is a
\emph{homology equivalence}, i.e.\ that the induced map in integral
homology is an isomorphism.  The Barratt-Priddy-Quillen theorem
(\cite{MR0314940}) gives a homology equivalence $\Z \times B
\Sigma_\infty \to QS^0$, where $QS^0$ is the infinite loop space
\begin{align*}
  QS^0 = \colim_{n \to \infty} \Omega^n S^n.
\end{align*}
The main theorem~\ref{thm:main} now takes the following equivalent
form.
\begin{theorem}\label{thm:MW}
  There is a homology equivalence
  \begin{align*}
    \Z \times B \Aut_\infty \to QS^0.
  \end{align*}
\end{theorem}

Alternatively the result can be phrased as a homotopy equivalence $\Z
\times B \Aut_\infty^+ \simeq QS^0$, where $B\Aut_\infty^+$ denotes
Quillen's plus-construction applied to $B\Aut_\infty$.  Quillen's
plus-construction converts homology equivalences to weak homotopy
equivalences, cf.\ e.g.\ \cite{MR649409}.

Most of the theorems stated or quoted above for $\Aut(F_n)$ have
analogues for mapping class groups.  Theorem~\ref{thm:HVW} above is
the analogue of the homological stability theorems of Harer and Ivanov
for the mapping class group (\cite{MR786348}, \cite{MR1015128}).
Corollary~\ref{cor:Mumford} above is the analogue of ``Mumford's
conjecture'', and the homotopy theoretic strengthening in
theorem~\ref{thm:MW} (which is equivalent to the statement in
theorem~\ref{thm:main}) is the analogue of Madsen-Weiss' generalized
Mumford conjecture (\cite{math.AT/0212321}, see also
\cite{math.AT/0605249}).

Some conjectures and partial results in this direction have been
known.  Hatcher (\cite{MR1314940}) noticed that there is a homotopy
equivalence $\Z \times B \Aut_\infty^+ \simeq QS^0 \times W$ for some
space $W$.  Hatcher and Vogtmann \cite{MR1671188} calculated
$H_k(\Aut(F_n);\Q)$ for small $k$.  They proved that
$H_4(\Aut(F_4);\Q) = \Q$ and that $H_k(\Aut(F_n);\Q) = 0$ for all
other $(k,n)$ with $0<k\leq 6$.  It follows that the stable rational
homology vanishes in degrees $\leq 6$, and they conjectured that
stable rational cohomology vanishes in all degrees.
Corollary~\ref{cor:Mumford} verifies Hatcher-Vogtmann's conjecture.

\subsection{Outline of proof}
\label{subsec:outline-proof}

Culler-Vogtmann's \emph{Outer Space} plays the role for $\Out(F_n)$
that Teichm\"uller space plays for mapping class groups.  Since its
introduction in \cite{MR830040}, it has been of central importance in
the field, and firmly connects $\Out(F_n)$ to the study of
\emph{graphs}.  A point in outer space $X_n$ is given by a triple
$(G,g,h)$ where $G$ is a connected finite graph, $g$ is a metric on
$G$, i.e.\ a function from the set of edges to $[0,\infty)$ satisfying
that the sum of lengths of edges in any cycle of $G$ is positive, and
$h$ is a marking, i.e.\ a conjugacy class of an isomorphism $\pi_1(G)
\to F_n$.  Two triples $(G,g,h)$ and $(G',g',h')$ define the same
point in $X_n$ if there is an isometry $\phi: G \to G'$ compatible
with $h$ and $h'$.  The isometry is allowed to collapse edges in $G$
of length 0 to vertices in $G'$.  If $G$ has $N$ edges, the space of
metrics on $G$ is an open subset $M(G)\subseteq [0,\infty)^N$.  Equip
$M(G)$ with the subspace topology and $X_n$ with the quotient topology
from $\amalg M(G)$, the disjoint union over all marked graphs $(G,h)$.
This defines a topology on $X_n$ and Culler-Vogtmann proves that it is
contractible.

Outer space is built using compact connected graphs $G$ with fixed
first Betti number $b_1(G) = n$.  The main new tool in this paper is
the definition of a space $\Phi(\R^N)$ of \emph{non-compact} graphs
$G\subseteq \R^N$.  Inside the space $\Phi(\R^N)$ is a space $B_N$ of
embedded compact graphs.  We will prove that a connected component of
$B_\infty$ is weakly equivalent to $B\Out(F_n)$.  Considering also
non-compact graphs allows us to define a map
\begin{align}\label{eq:2}
  B_N \xrightarrow{\tau_N} \Omega^N\Phi(\R^N).
\end{align}
In the analogy to mapping class groups, $\tau_N$ replaces the
Pontryagin-Thom collapse map of \cite{math.AT/0212321} and
\cite{math.AT/0605249} and as $N$ varies, the spaces $\Phi(\R^N)$ form
a spectrum $\mathbf{\Phi}$ which replaces the Thom spectrum $\MT{d}$
of \cite{math.AT/0605249} and $\CP^\infty_{-1}$ of
\cite{math.AT/0212321}.  We take the direct limit of~\eqref{eq:2} as
$N \to \infty$ and get a map
\begin{align*}
  B_\infty \xrightarrow{\tau_\infty} \Omega^\infty\mathbf{\Phi}
\end{align*}
or, by restriction to a connected component,
\begin{align*}
  B\Out(F_n) \to \Omega^\infty\mathbf{\Phi}.
\end{align*}
Composing with the map induced by the group maps
$\alpha_{n+1}^1\circ(\gamma_n^1\circ \beta_n^1): \Aut(F_n) \to
\Aut(F_{n+1}) \to \Out(F_{n+1})$ from~\eqref{eq:12} we get a map
\begin{align}\label{eq:3}
  \coprod_{n\geq 0} B \Aut(F_n) \to \Omega^\infty\mathbf{\Phi},
\end{align}
and the proof of theorem~\ref{thm:MW} is concluded in the following
steps.
\begin{enumerate}[(i)]
\item The map \eqref{eq:3} induces a well defined $\tau: \Z \times
  B\Aut_\infty \to \Omega^\infty\mathbf{\Phi}$,
\item $\tau$ is a homology equivalence,
\item $\Omega^\infty\mathbf{\Phi} \simeq QS^0$.
\end{enumerate}

The paper is organized as follows.  In chapter \ref{sec:sheaf-graphs}
we define and study the space $\Phi(\R^N)$.  In
chapter~\ref{sec:homotopy-types-graph} we define a subspace $B_N
\subseteq \Phi(\R^N)$ consisting of compact graphs and explain its
relation to $B\Out(F_n)$.  We also define the map~\eqref{eq:2}.  The
proof that there is an induced map $\tau: \Z \times B\Aut_\infty \to
\Omega^\infty\mathbf{\Phi}$ which is a homology equivalence is in
chapter~\ref{sec:graph-cobordism-category} and is in two steps.
First, in section~\ref{subsec:poset-model-graph} we define a
topological category $\mathcal{C}$, whose objects are finite sets and
whose morphisms are certain \emph{graph cobordisms}.  We prove the
equivalence $\Omega B \mathcal{C} \simeq \Omega^\infty \mathbf{\Phi}$.
Secondly, in section~\ref{subsec:posit-bound-subc}, we prove that
there is a homology equivalence $\Z \times B \Aut_\infty \to \Omega B
\mathcal{C}$.  This is very similar to, and inspired by, the
corresponding statements for mapping class groups in
\cite{math.AT/0605249}.  Finally, chapter~\ref{sec:homotopy-type-Phi}
is devoted to proving that $\Omega^\infty \mathbf{\Phi} \simeq QS^0$.
This completes the proof of theorem~\ref{thm:MW}.

In the supplementary chapter~\ref{sec:some-remarks-manif} we compare
with the work in \cite{math.AT/0605249}.  Our proof of
theorem~\ref{thm:MW} works with minor modifications if the space
$\Phi(\R^N)$ is replaced throughout by a space $\Psi_d(\R^N)$ of
smooth $d$-manifolds $M\subseteq \R^N$ which are closed subsets.  In
that case we prove an unstable version of the main result of
\cite{math.AT/0605249}.  To explain it, let $\Gr_d(\R^N)$ be the
Grassmannian of $d$-planes in $\R^N$, and $U^\perp_{d,N}$ the
canonical $(N-d)$ dimensional vector bundle over it.  Let
$\thh(U_{d,N}^\perp)$ be its Thom space.  Then we prove the weak
equivalence
\begin{align}
  \label{eq:26}
  B \mathcal{C}_d^N \simeq \Omega^{N-1} \thh(U_{d,N}^\perp),
\end{align}
where $\mathcal{C}_d^N$ is now the cobordism category whose objects
are closed $(d-1)$-manifolds $M \subseteq \{a\} \times \R^{N-1}$ and
whose morphisms are compact $d$-manifolds $W \subseteq [a_0, a_1]
\times \R^{N-1}$, cf.\ \cite[section 2]{math.AT/0605249}.  In the
limit $N \to \infty$ we recover the main theorem of
\cite{math.AT/0605249}, but~\eqref{eq:26} holds also for finite $N$.

\textbf{Acknowledgements.}  I am grateful to Ib Madsen for many
valuable discussions and comments throughout this project and for
introducing me to this problem as a graduate student in Aarhus; and to
Allen Hatcher, Kiyoshi Igusa and Anssi Lahtinen for useful comments on
earlier versions of the paper.

\section{The sheaf of graphs}
\label{sec:sheaf-graphs}

This chapter defines and studies a certain sheaf $\Phi$ on $\R^N$.
Roughly speaking, $\Phi(U)$ will be the set of all graphs $G \subseteq
U$.  We allow non-compact, and possibly infinite, graphs.  The precise
definition is given in section~\ref{subsec:sheaf-of-graphs} below,
where we also define a topology on $\Phi(U)$, making $\Phi$ a sheaf of
topological spaces.

\subsection{Definitions}
\label{subsec:sheaf-of-graphs}

Recall that a continuous map $f: X \to Y$ is a \emph{topological
  embedding} if $X \to f(X)$ is a homeomorphism, when $f(X)\subseteq
Y$ has the subspace topology.  If $X$ an $Y$ are smooth manifolds,
then $f$ is a $C^1$ embedding if $f$ is $C^1$, if $Df(x): T_xX
\to T_{f(x)}Y$ is injective for all $x \in X$, and if $f$ is a
topological embedding.

\begin{definition}\label{definition:Phi-U-as-set}
  Let $U\subseteq \R^N$ be open.  Let $\Phi(U)$ be the set of pairs
  $(G,l)$, where $G \subseteq U$ is a subset and $l: G \to [0,1]$ is a
  continuous function.  $(G,l)$ is required to satisfy the following
  properties.
  \begin{enumerate}[(i)]
  \item\label{item:7} If $p \in U - G$, there is a neighborhood $U_p
    \subseteq U$ of $p$ with $G \cap U_p = \emptyset$.
  \item\label{item:8} If $p \in l^{-1}([0,1))$, there is a
    neighborhood $U_p \subseteq U$ of $p$, an open set $V\subseteq
    (-1,1)$, and an $C^1$ embedding $\gamma: V \to U_p$ such that $G
    \cap U_p = \gamma(V)$ and $l\circ\gamma(t) = t^2$.
  \item\label{item:9} If $p \in l^{-1}(1)$, there is a neighborhood
    $U_p\subseteq U$ of $p$, a natural number $n \geq 3$, an open set
    $V \subseteq \vee^n [-1,1)$, and a topological embedding $\gamma:
    V \to U_p$ such that $G \cap U_p = \gamma(V)$ and $l\circ\gamma(t)
    = t^2$.  We require $\gamma$ to be $C^1$ in the following sense.
    Let $j_k: [-1,1) \to \vee^n [-1,1)$ be the inclusion of the $k$th
    wedge summand.  Then $\gamma_k = \gamma \circ j_k: (j_k)^{-1}(V_p)
    \to U_p$ is $C^1$, and the vectors $\gamma_k'(0)/|\gamma_k'(0)|$
    are pairwise different, $k = 1,\dots, n$.
  \end{enumerate}
  A \emph{graph in $U$} is an element $G \in \Phi(U)$.
\end{definition}
(\ref{item:7}) is equivalent to $G\subseteq U$ being a closed subset.
$G$ need not be compact, and non-compact $G \in \Phi(U)$ may have
infinitely many edges and vertices.  $G$ also need not be connected.

Let us say that a parametrizations $\gamma: V \to G$ as in
(\ref{item:8}), with $V\subseteq (-1,1)$ and satisfying $l\circ
\gamma(t) = t^2$, is an \emph{admissible parametrization of $G$ at
  $p$}.  These are almost unique: If $\bar \gamma$ is another
admissible parametrization at $p$, then either $\gamma(t) = \bar
\gamma(t)$ or $\gamma(t) = \bar \gamma(-t)$ for $t$ near
$\gamma^{-1}(p)$.  So specifying the function $l: G \to [0,1]$ is
equivalent to specifying an equivalence class $\{[\gamma],
[\bar\gamma]\}$ of germs of parametrizations around each point.  Often
we will omit $l$ from the notation and write e.g.\ $G \in \Phi(U)$.

An inclusion $U \subseteq U'$ induces a restriction map $\Phi(U') \to
\Phi(U)$ given by
\begin{align*}
  G \mapsto G \cap U.
\end{align*}
This makes $\Phi$ a sheaf on $\R^N$.  More generally, if $j:
U \to U'$ is an embedding (not necessarily an inclusion) of open
subsets of $\R^N$, define $j^*: \Phi(U') \to \Phi(U)$ by
\begin{align*}
  j^*(G) = j^{-1}(G)
\end{align*}
and $j^*(l) = l \circ j: j^*(G) \to [0,1]$.

We have the following standard terminology.
\begin{definition}
  Let $G \in \Phi(U)$.
  \begin{enumerate}[(i)]
  \item Let $\VV(G)$ = $l^{-1}(1)$.  This is the set of
    \emph{vertices} of $G$.
  \item An \emph{edge point} is a point in the 1-manifold $\EE(G) = G
    - \VV(G)$.
  \item An \emph{oriented edge} is a continuous map $\gamma: [-1,1]
    \to G$ such that $l\circ \gamma(t) = t^2$ and such that
    $\gamma|(-1,1)$ is an embedding.
  \item A \emph{closed edge} of $G$ is a subset $I\subseteq G$ which
    is the image of some oriented edge.  Each edge is the image of
    precisely two oriented edges.  If $I$ is the image of $\gamma$, it
    is also the image of the oriented edge given by $\bar \gamma(t) =
    \gamma(-t)$.
  \item A subset $T\subseteq G$ is a \emph{tree} if it is the union of
    finitely many vertices and closed edges of $G$ and if $T$ is
    contractible.
  \end{enumerate}
\end{definition}

The following notion of maps between elements of $\Phi(U)$ is
important for defining the topology on $\Phi(U)$.  We remark that it
does not make $\Phi(U)$ into a category (because composition is only
partially defined).
\begin{definition}\label{defn:morphisms}
  Let $G', G \in \Phi(U)$.  A morphism $\phi: G' \dashrightarrow G$ is
  a triple $(V',V,\phi)$, where $V'\subseteq G'$ and $V\subseteq G$
  are open subsets and $\phi: V' \to V$ is a continuous surjection
  satisfying the following two conditions.
  \begin{enumerate}[(i)]
  \item For each $v \in \VV(G)\cap V$, $\phi^{-1}(v)\subseteq G'$ is a
    tree.  Let $\VV(V) = \VV(G) \cap V$, $\EE(V) = V - \VV(V)$,
    \begin{align*}
      \VV(\phi) = \bigcup_{v\in \VV(V)} \phi^{-1}(v)
    \end{align*}
    and $\EE(\phi) = V' - \VV(\phi)$.
  \item $\phi$ restricts to a $C^1$ diffeomorphism of manifolds over
    $[0,1)$
    \begin{align}\label{eq:46}
      \EE(\phi) \to \EE(V).
    \end{align}
  \end{enumerate}
\end{definition}

Throughout the paper we will use dashed arrows for partially defined
maps.  Thus the notation $f: X \dashrightarrow  Y$ means that
$f$ is a function $f: U \to X$ for some subset $U\subseteq X$.

It can be seen in the following way that for any morphism
$(V',V,\phi)$ as above, the underlying map of spaces $\phi: V' \to V$
is proper.  Let $K\subseteq V$ be compact, and let $x_n \in
\phi^{-1}(K)$ be a sequence of points.  After passing to a subsequence
we can assume that the sequence $y_n = \phi(x_n)$ converges to a point
$y \in K$.  If $y \in \EE(V)$ we must have $x_n \to x$ with $x =
\phi^{-1}(y) \in \phi^{-1}(K)$ because~\eqref{eq:46} is a
diffeomorphism.  If $y \in \VV(V)$ then $T = \phi^{-1}(y)\subseteq G$
is a tree and it is immediate from the definitions that a tree has a
compact neighborhood $C\subseteq V'$ with $C - T \subseteq \EE(\phi)$.
It follows that $C = \phi^{-1}(\phi(C))$ and that $\phi(C)\subseteq V$
is a neighborhood of $V$.  Therefore $x_n \in C$ eventually, and hence
$x_n$ has a convergent subsequence.

\begin{definition}\label{defn:topology}
  Let $\epsilon > 0$.  Let $K\subseteq U$ be compact and write
  $K^\epsilon$ for the set of $k \in K$ with $\dist(k,U-K) \geq
  \epsilon$.
  \begin{enumerate}[(i)]
  \item $\phi = (V',V,\phi)$ is $(\epsilon,K)$-small if $K\cap
    G\subseteq V$, $K^\epsilon \cap G' \subseteq \phi^{-1}(K)$, and if
    \begin{align*}
      |k - \phi(k)| < \epsilon\quad \text{for all $k \in
        \phi^{-1}(K)$}.
    \end{align*}
    We point out that $\phi^{-1}(K)\subseteq V'$ is a compact set
    containing $K^\epsilon \cap G'$.
  \item If $Q\subseteq K - \VV(G)$ is compact, then $\phi$ is
    $(\epsilon,K,Q)$-small if it is $(\epsilon,K)$-small and if for
    all $q \in Q \cap G$ there is an admissible parametrization
    $\gamma$ with $q = \gamma(t)$ and
    \begin{align*}
      |(\phi^{-1}\circ \gamma)'(t) - \gamma'(t)| < \epsilon.
    \end{align*}
  \item For $\epsilon,K,Q$ as above, let $\UU_{\epsilon,K,Q}(G)$ be
    the set
    \begin{align*}
      \{G' \in \Phi(U) \mid \text{there exists an
        $(\epsilon,K,Q)$-small $\phi: G' \dashrightarrow G$}
      \}.
    \end{align*}
    For the case $Q = \emptyset$ we write $\UU_{\epsilon,K}(G) =
    \UU_{\epsilon,K,\emptyset}(G)$.
  \item The $C^0$-topology on $\Phi(U)$ is the topology generated by
    the set
    \begin{align}\label{item:15}
      \{\UU_{\epsilon,K}(G) \mid & G \in \Phi(U), \text{$\epsilon >
        0$, $K\subseteq U$ compact}\}.
    \end{align}
  \item The $C^1$-topology on $\Phi(U)$ is the topology generated by
    the set
    \begin{align}\label{item:14}
      \{\UU_{\epsilon,K,Q}(G) \mid & G \in \Phi(U), \text{$\epsilon >
        0$, $K\subseteq U$ and $Q \subseteq K - \VV(G)$ compact}\}.
    \end{align}
  \end{enumerate}
\end{definition}

Unless explicitly stated otherwise, we topologize $\Phi(U)$ using the
$C^1$-topology.  In lemma~\ref{lemma:basis} below we prove that the
sets~\eqref{item:15} and \eqref{item:14} form bases for the topologies
they generate, and that the sets $\UU_{\epsilon,K}(G)$ form a
neighborhood basis at $G$ in the $C^0$-topology for fixed $G$ and
varying $\epsilon > 0$, $K\subseteq U$ and similarly
$\UU_{\epsilon,K,Q}(G)$ in the $C^1$ topology.

\begin{example}\label{example:convergence-to-empty}
  We discuss the important example $G = \emptyset \in \Phi(U)$, which
  will illustrate the role of the compact set $K$.  Any morphism
  $(V',V,\phi): G' \dashrightarrow {\emptyset}$ must have $V' = V =
  \emptyset$, because $V\subseteq G = \emptyset$ and $\phi: V' \to V$.
  Thus $\phi$ is $(\epsilon,K)$-small if and only if $K^\epsilon\cap
  G' = \emptyset$.  In particular
  \begin{itemize}
  \item If $X$ is a topological space then a map $f: X \to \Phi(U)$ is
    continuous at a point $x \in X$ with $f(x) = \emptyset$ if and
    only if for all compact subsets $K \subseteq U$ there exists a
    neighborhood $V \subseteq U$ of $x$ such that $f(y) \cap K =
    \emptyset$ for all $y \in V$.
  \item If $(G_n)_{n \in \N}$ is a sequence of elements of $\Phi(U)$,
    then $G_n \to \emptyset$ if and only if for all compact subsets $K
    \subseteq U$ there exists $N \in \N$ such that $G_n \cap K =
    \emptyset$ for $n > N$.
  \end{itemize}
\end{example}

\begin{lemma}\label{lemma:Phi-is-connected}
  The space $\Phi(\R^N)$ is path connected.
\end{lemma}
\begin{proof}
  We construct an explicit path from a given $G \in \Phi(\R^N)$ to the
  basepoint $\emptyset \in \Phi(\R^N)$.  Choose a point $p \in \R^N -
  G$ and let $\phi_t: \R^N \to \R^N$, $t\in [0,1]$ be the map given by
  \begin{align*}
    \phi_t(x) = (1-t)x + tp.
  \end{align*}
  Then $\phi_t$ is a diffeomorphism for $t < 1$ and $\phi_1(x) = p$
  for all $x$.  Let $G_t = (\phi_t)^{-1}(G)$.  This defines a map $t
  \mapsto G_t \in \Phi(\R^N)$.  We will see later
  (lemma~\ref{prop:res-is-cont}) that it is continuous on $[0,1)$.
  Continuity at 1 can be seen as follows.  For a given compact
  $K\subseteq \R^N$, choose $\delta > 0$ such that $K \subseteq
  B(p,\delta^{-1})$.  Then $G_t \cap K = \emptyset$ for $t >
  1-\delta$.
\end{proof}

\subsection{Point-set topological properties}
\label{sec:point-set-topol}

In this section we prove various results about $\Phi(U)$ of a
point-set topological nature.  The verifications are elementary, but
somewhat tedious, and their proofs could perhaps be skipped at a first
reading.

Let us first point out that $V$ in definition~\ref{defn:topology} can
always be made smaller: If $(V',V,\phi): G' \dashrightarrow G$ is
$(\epsilon,K,Q)$-small, then $(W',W,\psi)$ is $(\epsilon,K,Q)$-small
if $K\cap G \subseteq W \subseteq V$ with $W\subseteq V$ open, $W' =
\phi^{-1}(W)$, and $\psi = \phi|W'$.

\begin{lemma}
  \label{lemma:triangle}
  If $(V',V,\phi): G' \dashrightarrow G$ is $(\epsilon,K,Q)$-small and
  $(V'',W',\psi): G'' \to G'$ is $(\delta,K',Q')$-small with $K'
  \supseteq \phi^{-1}(K) \cup K^\epsilon$ and $Q' \supseteq
  \phi^{-1}(Q)$, then $(V'',V,\phi\circ \psi)$ is an
  $(\epsilon+\delta,K)$-small morphism after possibly shrinking $V$
  and $W'$.
\end{lemma}
\begin{proof}
  It suffices to consider the case $K' = \phi^{-1}K \cup K^\epsilon$.
  We have $\phi^{-1}(K) \subseteq K' \cap G'\subseteq W'$ by
  assumption on $K'$ and by $(\delta,K')$-smallness of
  $(V'',W',\psi)$.  
  Therefore the subset $\phi(V'-W') \subseteq V$ is disjoint from $K$,
  and properness of $\phi: V' \to V$ implies that $\phi(V'-W')
  \subseteq V$ is closed.  After replacing $V$ by $V - \phi(V'-W')$ we
  can assume that $V' = \phi^{-1}(V) \subseteq W'$.

  We have $K^\epsilon \cap G'\subseteq \phi^{-1}(K) \subseteq V'$
  because $(V',V,\phi)$ is $(\epsilon,K)$-small, so $K' \cap G'
  \subseteq (K^\epsilon \cup \phi^{-1}(K))\cap G' \subseteq V'$.
  Hence after shrinking $W'$ we can assume that $W' = V'$.  Then
  $(V'',V,\phi\circ \psi)$ is a morphism.  It is $(\epsilon +
  \delta,K,Q)$-small because $K^{\epsilon + \delta} \cap G'' \subseteq
  (\phi \circ \psi)^{-1}(K)$ and for $k \in (\phi\circ \psi)^{-1}(K)
  \subseteq \psi^{-1}(K')$ we have
  \begin{align*}
    |k - \phi\circ \psi(k)| \leq |k - \psi(k)| + |\psi(k) -
    \phi(\psi(k))| < \delta + \epsilon.
  \end{align*}
  A similar condition on first derivatives holds on $(\phi\circ
  \psi)^{-1}(Q)$.
\end{proof}

\begin{lemma}\label{lemma:basis}
  Let $G \in \Phi(U)$, $\epsilon > 0$ and $K\subseteq U$ and $Q
  \subseteq K - \VV(G)$ compact.  Let $G' \in \UU_{\epsilon,K,Q}(G)$.
  Then there exists $\delta > 0$ and compact $K'\subseteq U$,
  $Q'\subseteq K' - \VV(G'$) such that
  \begin{align}\label{eq:50}
    \UU_{\delta,K',Q'}(G') \subseteq \UU_{\epsilon,K,Q}(G).
  \end{align}
  We can take $Q' = \emptyset$ if $Q = \emptyset$.
\end{lemma}
\begin{proof}
  Let $(V', V, \phi): G' \dashrightarrow G$ be $(\epsilon,K,Q)$-small.
  By compactness of $\phi^{-1}(K)$, we can choose $\delta > 0$ with
  $|\phi(k) - k| < \epsilon - \delta$ for all $k \in \phi^{-1}(K)$.
  By compactness of $Q$ we can assume that
  \begin{align*}
    |(\phi^{-1}\circ \gamma)'(t) - \gamma'(t)| < \epsilon - \delta.
  \end{align*}
  for all admissible parametrizations $\gamma$ of $G$ with $\gamma(t)
  = q \in Q$.  We can also assume that $\delta$ satisfies $\delta <
  \dist(K^\epsilon, \phi^{-1}(G - \Int(K))$.  Then $(V',V,\phi)$ is
  actually $(\epsilon-\delta, K, Q)$-small, and the claim follows from
  lemma~\ref{lemma:triangle} if we set
  \begin{align}
    K' = \phi^{-1}(K) \cup K^{(\epsilon-\delta)}, \quad Q' =
    \phi^{-1}(Q). \tag*{\qedhere}
  \end{align}
\end{proof}
Lemma~\ref{lemma:basis} implies that the set \eqref{item:14} is a
basis for the topology it generates, and that the collection of
$\UU_{\epsilon,K,Q}(G)$ forms a neighborhood basis at $G$, for fixed
$G$ and varying $\epsilon, K, Q$.  Similarly for the $C^0$-topology.

The next lemma is the main rationale for including the map $l: G \to
[0,1]$ into the data of an element of $\Phi(U)$.  It gives a partial
uniqueness result for the $(\epsilon,K)$-small maps $(V',V,\phi):
G' \dashrightarrow G$, whose existence is assumed when $G'
\in \UU_{\epsilon,K}(G)$.
\begin{lemma}\label{lemma:painful-proof}
  For each $G \in \Phi(U)$ and each compact $C \subseteq U$, there
  exists an $\epsilon > 0$ and a compact $K \subseteq U$ with
  $C\subseteq \Int(K^\epsilon)$ such that for $G' \in
  \UU_{\epsilon,K}(G)$, any two $(\epsilon,K)$-small maps
  \begin{align*}
    \phi, \psi: G' \dashrightarrow G
  \end{align*}
  must have $\phi = \psi$ near $C$.
\end{lemma}

Notice that both $\psi$ and $\phi$ are defined near $C \cap G'$ when
$C \subseteq \Int(K^\epsilon)$.

\begin{proof}
  Let $W\subseteq U$ be an open set with compact closure $\overline W
  \subseteq U$ and $C \subseteq W$.  We will prove that $\epsilon$ and
  $K$ can be chosen so that $\phi^{-1}(w) = \psi^{-1}(w)$ for all $w
  \in W$ when both maps are $(\epsilon,K)$-small.  If we also arrange
  $\epsilon < \dist(C,U-W)$, that will prove the statement in the
  lemma.

  First take $\epsilon > 0$ and $K\subseteq U$ such that $W \subseteq
  K^{2\epsilon}$.  We can assume that the distance between any two
  elements of $K \cap \VV(G)$ is greater than $2\epsilon$.  Then the
  triangle inequality implies that $\phi^{-1}(v) = \psi^{-1}(v)$ for
  all $v \in K^{2\epsilon} \cap \VV(G)$.  It remains to treat edge
  points.

  Let $M \subseteq G \cap \Int(K)$ be the smallest open and closed
  subset containing $G \cap W$.  We claim that $\phi^{-1}(v)=
  \psi^{-1}(v)$ for $v$ in $M - \VV(G)$.  It suffices to consider $v
  \in M \cap l^{-1}((0,1))$ since that set is dense in $M - \VV(G)$
  (we omit only ``midpoints'' of edges).  Compactness of $\overline W$
  implies that $\pi_0 M$ is finite (connected components of $G \cap
  \Int(K)$ are open in $G \cap \Int(K)$, so the compact subset
  $\overline W \cap G$ can be non-disjoint from only finitely many).
  $l^{-1}(\{0,1\}) \cap K$ is a finite set of points, so $\pi_0(M \cap
  l^{-1}((0,1)))$ is also finite.  Choose a $\tau > 0$ such that the
  inclusion
  \begin{align}\label{eq:51}
    W \cap l^{-1}([\tau,1-\tau]) \to M \cap l^{-1}((0,1))
  \end{align}
  is a $\pi_0$-surjection (i.e. the induced map on $\pi_0$ is
  surjective).

  The function $l': G' \to [0,1]$ restricts to a local diffeomorphism
  $(l')^{-1}((0,1)) \to (0,1)$.  It follows that the diagonal embedding
  \begin{align}
    \label{eq:52}
    (l')^{-1}((0,1)) \xrightarrow{\mathrm{diag}} \{(k,m) \in G'\times G'
    \mid l(k) = l(m)
    \in (0,1)\}
  \end{align}
  has open image.  Therefore (by continuity of $\psi^{-1}, \phi^{-1}:
  l^{-1}((0,1)) \to G'$) the set
  \begin{align*}
    \{k \in l^{-1}((0,1)) \cap M \mid \psi^{-1}(k) = \phi^{-1}(k) \}
  \end{align*}
  is open and closed in $l^{-1}((0,1)) \cap M$, so it suffices to
  prove that it contains a point in each path component of
  $l^{-1}((0,1))\cap M$.  We prove that $\psi^{-1} = \phi^{-1}$ when
  composed with the $\pi_0$-surjection \eqref{eq:51}.

  The set of pairs $(k,m)$ with $k,m \in K\cap l^{-1}([\tau,1-\tau])$,
  and $l(k) = l(m)$ and $k \neq m$ is a compact subset of $K \times
  K$, so we can assume that $|k-m| > 2\epsilon$ for such $(k,m)$.  Now
  let $k \in W \cap l^{-1}([\tau,1-\tau]) \in K^{2\epsilon} \cap G'$,
  assume $\psi^{-1}(k) \neq \phi^{-1}(k)$, and set $x = \psi^{-1}(k)$.
  $\psi$ is $(\epsilon,K)$-small, so $|x-k| = |x - \psi(x)| <
  \epsilon$.  Hence $x \in K^\epsilon\cap G'$, so $\phi(x)$ is defined
  and $\phi(x) \in G \cap G$.  Injectivity of $\phi$ (on non-collapsed
  edges) implies that $\phi(x) = \phi(\psi^{-1}(k)) \neq k = \psi(x)$.
  Set $m = \phi(x)$.  Since $l(k) = l'(x) = l(m)\in [\tau, 1-\tau]$,
  we have
  \begin{align*}
    2\epsilon < |k-m| = |\psi(x) - \phi(x)| \leq |\psi(x) - x| + |x -
    \phi(x)|
  \end{align*}
  which contradicts $\phi$ and $\psi$ being $(\epsilon,K)$-small.
\end{proof}

\begin{proposition}\label{prop:sheaf-of-spaces}
  $\Phi$ is a sheaf of topological spaces on $\R^N$, i.e. the
  following diagram is an equalizer diagram of topological spaces for
  each covering of $U$ by open sets $U_j$, $j \in J$
  \begin{align*}
    \Phi(U) \to \prod_{j\in J} \Phi(U_j) \rightrightarrows
    \prod_{(i,l)\in J\times J} \Phi(U_i \cap U_l).
  \end{align*}
\end{proposition}
\begin{proof}
  Let $V\subseteq U$ be open and let $r: \Phi(U) \to \Phi(V)$ denote
  the restriction map.  If $G \in \Phi(V)$ and $G' \in
  r^{-1}\UU_{\epsilon,K,Q}(G)$, lemma~\ref{lemma:basis} provides
  $\delta > 0$ and $K',Q'\subseteq V$ such that
  $\UU_{\delta,K',Q'}(rG') \subseteq \UU_{\epsilon,K,Q}(G)$.  If
  $\dist(K, \R^N - V) > \delta$ we have
  \begin{align*}
    \UU_{\delta,K',Q'}(G') = r^{-1}\UU_{\delta,K',Q'}(rG') \subseteq
    r^{-1}\UU_{\epsilon,K,Q}(G)
  \end{align*}
  which proves that $r^{-1}\UU_{\epsilon,K,Q}(G)$ is open and hence
  that $r$ is continuous.  Therefore the maps in the diagram are all
  continuous.  The proposition holds for both the $C^0$ and the $C^1$
  topology.  We treat the $C^0$ case first.

  Let $\tilde \Phi(U)$ denote the image of $\Phi(U) \to \prod
  \Phi(U_i)$, topologized as a subspace of the product.  Then $\Phi(U)
  \to \tilde \Phi(U)$ is a continuous bijection.  Take $G \in \Phi(U)$
  and $\epsilon > 0$ and let $K\subseteq U$ be compact.  We will prove
  that the image of $\UU_{\epsilon,K}(G)\subseteq \Phi(U)$ in $\tilde
  \Phi(U)$ is a neighborhood of the image $\tilde G\in \tilde \Phi(U)$
  of $G\in \Phi(U)$.

  Choose a finite subset $\{j_1, \dots, j_n\} \subseteq J$ and compact
  $C_i\subseteq U_{j_i}$ such that $K\subseteq \cup_{i=1}^n
  C_i^\delta$ for some $\delta \in (0,\epsilon)$.  Let $K_i \subseteq
  U_{j_i}$ be compact subsets with $C_i \subseteq \Int(K_i)$.  Let
  $K_{il} = K_i \cap K_l$.  By lemma~\ref{lemma:painful-proof} we can
  assume, after possibly shrinking $\delta$ and enlarging the $K_i$,
  that $(\delta,K_{il})$-small morphisms $\phi_{il}: G_{il}
  \dashrightarrow (G|U_{j_ij_l})$ with $G_{il} \in \Phi(U_{il})$ have
  unique restriction to a neighborhood of $G \cap C_{il}$.  Thus, if
  $G' \in \Phi(U)$ has a $(\delta,K_i)$-small $\phi_{i}: (G'|U_{j_i})
  \dashrightarrow (G|U_{j_i})$ for all $i = 1,\dots, n$, then $\phi_i$
  and $\phi_l$ agrees near $G \cap C_{il}$.  Therefore they glue to a
  morphism $\phi: G' \dashrightarrow G$ which is defined near $L =
  \cup_i C_i$ and agrees with $\phi_i$ near $C_i$.  Since $\phi_i$ is
  $(\delta,K_i)$-small we will have $\phi_i(C_i) \supseteq C_i^\delta
  \cap G$ and hence $\phi(L) \supseteq K \cap G$ so the image of
  $\phi$ contains $K \cap G$.  The domain contains $L \cap G' =
  \cup_{i=1}^n (C_i \cap G')$ which contains $G \cap G'$ and hence
  $K^\delta \cap G'$.  Finally, let $k \in G$ have $\phi(k) \in K$ and
  hence $\phi(k) C_i^\delta\subseteq K_i$ for some $i$.  Then $\phi(k)
  = \phi_i(k)$ and
  \begin{align*}
    |\phi(k) - k| = |\phi_i(k) - k| < \delta
  \end{align*}
  because $\phi_i$ is $(\delta,K_i)$-small.  We get that $\phi$ is
  $(\delta,K)$-small.

  We have proved that $G' \in \UU_{\delta,K}(G) \subseteq
  \UU_{\epsilon,K}(G)$ whenever $(G'|U_{j_i})\in
  \UU_{\delta,K_i}(G|U_{j_i})$ for each $i = 1, \dots, n$.  Therefore
  the image of $\UU_{\epsilon,K}(G)\subseteq \Phi(U)$ in $\tilde
  \Phi(U)$ contains $p^{-1}(\UU)$, where
  \begin{align}\label{eq:48}
    \UU = \prod_{i=1}^n \UU_{\delta,K_i}(G|U_{j_i}) \subseteq
    \prod_{i=1}^n \Phi(U_{j_i}),
  \end{align}
  and $p$ is the projection $p: \tilde \Phi(U) \to \prod_{i=1}^n
  \Phi(U_i)$.  $p^{-1}(\UU)$ is the required neighborhood of $\tilde
  G$.

  To prove the $C^1$ case, suppose $Q \subseteq K - \VV(G)$, repeat
  the proof of the $C^0$ case, and set $Q_i = K_i \cap Q$.  Then
  replace $\UU_{\delta,K_i}$ by $\UU_{\delta,K_i,Q_i}$
  in~\eqref{eq:48}.
\end{proof}

The sheaf property implies that continuity of a map $f: X \to \Phi(U)$
can be checked locally in $X\times U$.  In other words, $f$ is
continuous if for each $x \in X$ and $u \in U$ there is a neighborhood
$V_x \times W_u \subseteq X \times U$ such that the composition
\begin{align*}
  V_x \to X \xrightarrow{f} \Phi(U) \xrightarrow{\mathrm{restr.}}
  \Phi(W_u)
\end{align*}
is continuous.  In particular, $U \mapsto \Map(X, \Phi(U))$ is a sheaf
for every space $X$.

\begin{proposition}\label{prop:res-is-cont}
  If $V\subseteq U$ are open subsets of $\R^N$, then the restriction
  map $\Phi(U) \to \Phi(V)$ is continuous.  More generally, the action
  map
  \begin{align*}
    \Emb(V,U) \times \Phi(U) \to \Phi(V)
  \end{align*}
  $(j,G) \mapsto j^*(G)$ is continuous, where $\Emb(V,U)$ is given the
  $C^1$ topology.
\end{proposition}
\begin{proof}[Proof sketch.]
  Let $j \in \Emb(V,U)$, $G \in \Phi(U)$ and let $\epsilon > 0$, and
  let $K \subseteq V$ be compact.  Choose $\delta > 0$ and compact
  subsets $C\subseteq V$ and $L\subseteq jV$ such that $K \subseteq
  C^\delta$ and $jK \subseteq L^\delta$.  Choose a number $M$ such that
  \begin{align*}
    |j^{-1}(l) - j^{-1}(l')| \leq M|l - l'| \quad \text{and} \quad
    |D_lj^{-1}(v)| \leq M |v|
  \end{align*}
  for all $l,l' \in L$ and $v \in T_l(U)$.  We can assume $2M\delta
  \leq \epsilon$.  Then $j^*\phi: j^*G' \to j^*G$ is
  $(\epsilon/2,K)$-small if $\phi: G' \to G$ is $(\delta,K)$-small.

  Let $j' \in \Emb(V,U)$ be another embedding such that $j^{-1} \circ
  j': V \dashrightarrow V$ is $(\epsilon/2,K')$-small in the sense
  that the domain contains $(K')^{\epsilon/2}$, the image contains
  $K'$, and that $|f(k) - k| < \epsilon/2$ for $f(k) \in K'$.  Then
  the composition
  \begin{align*}
    (j')^*G' \xrightarrow{j^{-1}\circ j'} j^*(G')
    \xrightarrow{j^*\phi} j^*G
  \end{align*}
  is $(\epsilon,K)$-small by lemma~\ref{lemma:triangle}, provided $K'
  \supseteq (j^*\phi)^{-1}(K) \cup K^\epsilon$, which is satisfied if
  $(K')^\epsilon \supseteq K$.  This proves continuity when $\Phi(U)$
  and $\Phi(V)$ are given the $C^0$ topology.  The $C^1$ topology is
  similar.
\end{proof}

\section{Homotopy types of graph spaces}
\label{sec:homotopy-types-graph}

Lemma~\ref{lemma:Phi-is-connected} shows that the full space
$\Phi(\R^N)$ is path connected.  A similar argument shows that
$\Phi(\R^N)$ is in fact $(N-3)$-connected.  In this section we study
the homotopy types of certain subspaces of $\Phi(\R^N)$.

\subsection{Graphs in compact sets}
\label{sec:graphs-compact-sets}

\begin{definition}
  \begin{enumerate}[(i)]
  \item For a closed subset $A\subseteq U$, let $\Phi(A)$ be the set
    of germs around $A$, i.e.\ the colimit of $\Phi(V)$ over open sets
    with $A\subseteq V\subseteq U$.  We remark that the colimit
    topology is often not well behaved (for example if $A$ is a point
    then the one-point subset $\{[\emptyset]\} \subseteq \Phi(A)$ is
    dense), and we consider $\Phi(A)$ as a set only.
  \item Let $U\subseteq \R^N$ be open and $M\subseteq U$ compact.  For
    a germ $S \in \Phi(U-\Int M)$, let $\Phi^S(M)$ be the inverse
    image of $S$ under the restriction $\Phi(U) \to \Phi(U - \Int M)$.
    Topologize $\Phi^S(M)$ as a subspace of $\Phi(U)$.
  \item Let $G',G \in \Phi^S(M)$.  A \emph{graph epimorphism} $G' \to
    G$ is a morphism $(V',V,\phi)$ in the sense of
    definition~\ref{defn:morphisms} which is surjective and everywhere
    defined (i.e.\ $V' = G'$ and $V = G$).  Furthermore $\phi$ is
    required to restrict to the identity map $S \to S$.
  \item Let $\mathcal{G}_S$ be the category with objects $\Phi^S(M)$
    and graph epimorphisms as morphisms.  We consider
    $\Ob(\mathcal{G}_S)$ and $\mor(\mathcal{G}_S)$ discrete sets.
  \end{enumerate}
\end{definition}
The main result in this section is the following theorem.

\begin{theorem}\label{thm:Phi-is-B-Cat}
  Let $U\subseteq \R^N$ be open and $M\subseteq U$ compact.  Let $S
  \in \Phi(U-\Int M)$.  Assume $\Int M$ is $(N-3)$-connected.  Then
  there is an $(N-3)$-connected map
  \begin{align*}
    \Phi^S(M) \to B\mathcal{G}_S.
  \end{align*}
\end{theorem}

In the next section we prove that the classifying space
$B\mathcal{G}_S$ is homotopy equivalent to a space built out the
spaces $BA_n^s$, where $A_n^s$ are the groups from
theorem~\ref{thm:HVW}.  Combined with theorem~\ref{thm:Phi-is-B-Cat}
above this leads to theorem~\ref{thm:conclusion-htpy-types}, which
summarizes the results of sections~\ref{sec:graphs-compact-sets},
\ref{sec:spac-graph-embedd}, and \ref{sec:abstract-graphs}.  We need
more definitions for the proof.
\begin{definition}
  \begin{enumerate}[(i)]
  \item Let $M\subseteq U$ be compact and $R \in \Phi(U)$.  Let $S =
    [R] \in \Phi(U - \Int M)$ be the germ of $R$.  Let $\Phi(M;R)$ be
    the set of pairs $(G,f)$, where $G \in \Phi^S(M)$ and $f: G \to R$
    is a graph epimorphism.
  \item Let $(G,f) \in \Phi(M;R)$.  For $(\epsilon,K,Q)$ as in
    definition~\ref{defn:topology}, let
    $\UU_{\epsilon,K,Q}(G,f)\subseteq \Phi(M;R)$ be the set of
    $(G',f')$ which admits an $(\epsilon,K,Q)$-small $\phi: G'
    \dashrightarrow G$ with $f' = \phi \circ f$.  Topologize
    $\Phi(M;R)$ by declaring that the $\UU_{\epsilon,K,Q}(G,f)$ form a
    basis.
  \item Let $\Emb_R(M)\subseteq \Phi(M;R)$ be the subspace in which
    the morphism $f: G \to R$ has an inverse morphism $f^{-1}: R \to
    G$.
  \end{enumerate}
\end{definition}

The space $\Emb_R(M)$ can be thought of as a space of certain
embeddings $R \to U$.  Namely $(G,f)\in \Emb_R(M)$ can be identified
with the map $f^{-1}: R \to G \subseteq U$.

Throughout the paper we will make extensive use of simplicial spaces.
Recall that a simplicial space $X_\bullet$ has a geometric realization
$\|X_\bullet\|$ and that a simplicial map $f_\bullet: X_\bullet \to
Y_\bullet$ induces $\|f_\bullet\|: \|X_\bullet\| \to \|Y_\bullet\|$.
If for each $k$ the map $f_k: X_k \to Y_k$ is $(n-k)$-connected, the
geometric realization $\|f_\bullet\|$ is $n$-connected.  In particular
$\|f_\bullet\|$ is a weak equivalence if each $f_k$ is a weak
equivalence.  Recall also that to each category $C$ is an associated
classifying space $BC$, defined as the geometric realization of the
nerve $N_\bullet C$.

If $C$ is a category (usually not topologized) and $F: C \to
\mathrm{Spaces}$ is a functor, then the \emph{homotopy colimit} of $F$
is defined as
\begin{align*}
  \hocolim_C F = B(C \wr F)
\end{align*}
where $(C\wr F)$ is the category whose objects are pairs $(c,x)$ with
$c \in \Ob(C)$ and $x \in F(c)$, and whose morphisms $(c,x) \to
(c',x')$ is the set of morphisms $f \in C(c,c')$ with $F(f)(x) = x'$.
If $T: F \to G$ is a natural transformation such that $T(x): F(x) \to
G(x)$ is $n$-connected for each object $x$, the induced map $\hocolim
F \to \hocolim G$ is also $n$-connected.

The proof of theorem~\ref{thm:Phi-is-B-Cat} is broken down into the
following assertions, whose proofs occupy the remainder of this
section and the following.
\begin{itemize}
\item The forgetful map
  \begin{align*}
    \hocolim_{R \in \mathcal{G}_S} \Phi(M;R) \to \Phi^S(M)
  \end{align*}
  induced by the projection $(G,f) \mapsto G$ is a weak equivalence.
\item The inclusion $\Emb_R(M) \to \Phi(M;R)$ is a weak equivalence.
\item The space $\Emb_R(M)$ is $(N-4)$-connected if $\Int(M)$ is
  $(N-3)$-connected.
\end{itemize}

The following lemma will be used again throughout the paper.  Recall
that a map is \emph{etale} if it is a local homeomorphism and an open
map.
\begin{lemma}\label{lemma:Segal} 
  Let $C$ be a topological category and $Y$ a space.  Regard $Y$ as a
  category with only identity morphisms, and let $f: C \to Y$ be a
  functor such that $N_0 f$ and $N_1 f$ are etale maps.  Assume that
  $B(f^{-1}(y))$ is contractible for all $y \in Y$.  Then $Bf: BC \to
  Y$ is a weak equivalence.
\end{lemma}
\begin{proof}[Proof sketch]
  The hypothesis implies that a neighborhood of the fiber
  \begin{align*}
    B(f^{-1}(y)) \subseteq BC
  \end{align*}
  is homeomorphic, as a space over $Y$, to a neighborhood of
  \begin{align*}
    \{y\} \times B(f^{-1}(y)) \subseteq Y \times
    B(f^{-1}(y)).
  \end{align*}
  Then the result follows from \cite[proposition (A.1)]{MR516216}.
\end{proof}

\begin{lemma}\label{lemma:p-is-etale}
  The forgetful map $p: \Phi(M;R) \to \Phi^S(M)$, $p(G,f) = G$, is
  etale.
\end{lemma}
\begin{proof}
  Let $(G,f) \in \Phi(M;R)$.  An application of
  lemma~\ref{lemma:painful-proof} gives an $\epsilon > 0$ and a
  compact $K\subseteq U$ such that any $G' \in \Phi^S(M) \cap
  \UU_{\epsilon,K}(G)$ will have a \emph{unique} graph epimorphism
  $\phi_{G'}: G' \to G$ which is $(\epsilon,K)$-small.  $\phi_{G'}$
  restricts to the identity outside $M$.  This gives a map $G' \mapsto
  (G',f \circ \phi_{G'})$ which is a local inverse to $p$.  We have
  proved that $p$ restricts to a homeomorphism
  \begin{align}
    \UU_{\epsilon,K}(G,f) \to \UU_{\epsilon,K}(G).\tag*{\qedhere}
  \end{align}
\end{proof}

\begin{proposition}\label{prop:hocolim-weak-equiv}
  The map
  \begin{align*}
    \hocolim_{R \in \mathcal{G}_S} \Phi(M;R) \to \Phi^S(M)
  \end{align*}
  induced by the projection $(G,f) \mapsto G$ is a weak equivalence.
\end{proposition}
\begin{proof}
  The maps from lemma~\ref{lemma:p-is-etale} assemble to a map
  \begin{align*}
    \coprod_{R \in \mathcal{G}_S} \Phi(M;R) \xrightarrow{p}
     \Phi^S(M).
  \end{align*}
  The domain of this map is the space of objects of the category
  $(\mathcal{G}_S \wr \Phi(M;-))$.  Morphisms $(R',(G',f')) \to
  (R,(G,f))$ exist only if $G' = G$; then they are morphisms $\phi: R'
  \to R$ in $\mathcal{G}_S$ with $\phi \circ f' = f$.  The classifying
  space of this category is the homotopy colimit in the proposition,
  and $ p$ induces a map
  \begin{align*}
    B p: \hocolim_{R \in \mathcal{G}_S} \Phi(M;R) \to \Phi^S(M).
  \end{align*}
  Let $G \in \Phi^S(M)$.  Then the subcategory $ p^{-1}(G) \subseteq
  (\mathcal{G}_S \wr \Phi(M;-))$ has $(G,(G,\id))$ as initial object.
  Therefore $(B p)^{-1}(G) = B( p^{-1}(G))$ is contractible, so $ p$
  satisfies the hypotheses of lemma~\ref{lemma:Segal}.
\end{proof}

\subsection{Spaces of graph embeddings}
\label{sec:spac-graph-embedd}

\begin{figure}
  \centering
  \includegraphics{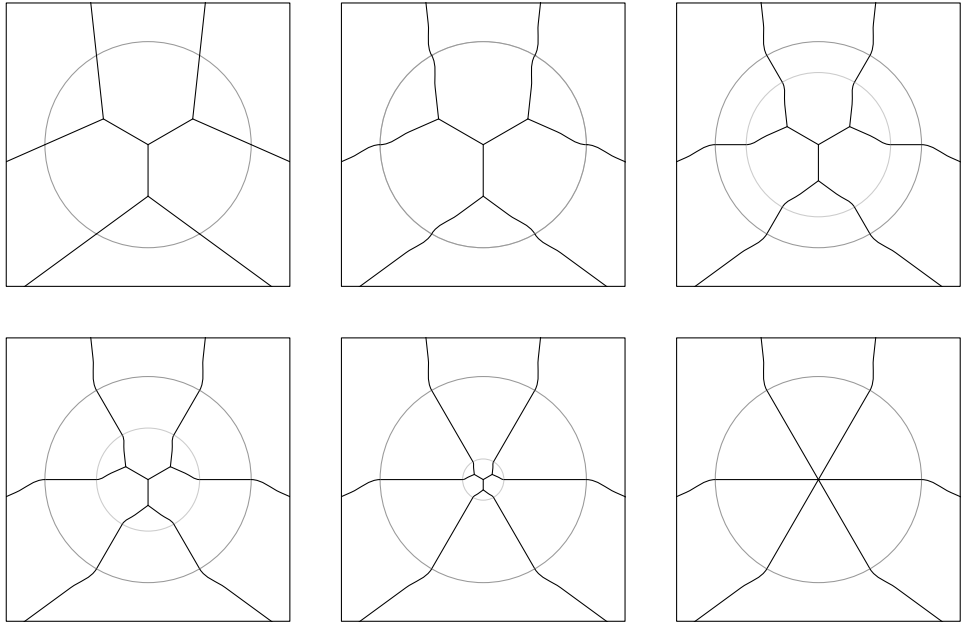}
  \caption{$G_t$ for various $t \in [0,1]$, cf.\
    construction~\ref{constr}.}
  \label{fig:2}
\end{figure}

Our next aim is to prove that the space $\Phi(M;R)$ is highly
connected when $\Int M \subseteq \R^N$ is highly connected and $N$ is
large.  The main step is to prove that the inclusion $\Emb_R(M) \to
\Phi(M;R)$ is a weak equivalence.  Although it is slightly lengthy to
give all details, the idea is easy to explain.  Suppose $(G,f) \in
\Phi(M;R)$, and we want to construct a path to an element in
$\Emb_R(M)$.  The map $f: G \to R$ specifies a finite set of trees
$T_v = f^{-1}(v)\subseteq G$, $v \in \VV(R) \cap \Int M$, such that
$G$ becomes isomorphic to $R$ when every $T_v\subseteq G$ is collapsed
to a point.  The point is that this contraction can be carried out
inside $M$, by continuously shortening leaves of the tree $T_v$ and
``dragging along'' edges incident to $T_v$ (see the illustration in
figure~\ref{fig:2}).  The formal proof consists of making this
construction precise and proving it can be done continuously.  The
construction is remotely similar to the Alexander trick.

We begin by constructing a prototype collapse.  This is done in
construction~\ref{constr} below, illustrated in figure~\ref{fig:2}.
\begin{definition}\label{definition:collapsible-pos}
  Let $(G,l) \subseteq \Phi(\R^N)$, and let $T \subseteq G$ be a tree.
  An \emph{incident edge} to $T$ is a map $\gamma:[0,\tau]\to G$ with
  $\tau < 2$ such that $l(\gamma(t)) = (t-1)^2$ and $\gamma^{-1}(T) =
  \{0\}$.  We consider two incident edges equivalent if one is a
  restriction of the other.  Say that $(G,T)$ is in \emph{collapsible
    position} if all $g \in G \cap B(0,3)$ are in either $T$ or in the
  image of an incident edge, if $T \subseteq \Int D^N$, and if there
  are representatives $\gamma_i:[0,\tau_i] \to G$ for all the incident
  edges satisfying
  \begin{align*}
    |\gamma_i(\tau_i)| & > 3,\\
    \langle \gamma_i(t), \gamma_i'(t) \rangle & \geq 0, \quad \text{when
      $|\gamma_i(t)| \in [1, 3]$}.
  \end{align*}
  These $\gamma_i$ provide a ``distance to $T$'' function $d: G \cap
  B(0,3) \to [0,2)$ given by $d(x) = 0$ when $x \in T$ and
  $d(\gamma_i(t)) = t$.
\end{definition}
We point out that if $G \in \Phi(\R^N)$ and if there exists a $T$ with
$(G,T)$ in collapsible position, then $T$ is unique (it must be the
union of all closed edges of $G$ contained in $\Int B(0,1)$), and the
function $d: G \cap B(0,3) \to [0,2)$ is independent of choice of
representatives $\gamma_i$.

\begin{construction}\label{constr}
  Let $\lambda_{\frac13}: [0,\infty) \to [0,\infty)$ be a smooth
  function satisfying $\lambda_{\frac13}(r) = 3r/2$ for $r \leq 1.3$,
  $\lambda_{\frac13}(r) = 2$ for $1.4 \leq r \leq 1.9$, and
  $\lambda_{\frac13}(r) = r$ for $r > 2.5$.  We also assume
  $\lambda_\frac13'(r) \geq 0$ and $\lambda_\frac13'(r) > 0$ for
  $\lambda_\frac13(r) \neq 2$ and $\lambda_\frac13'(r)\leq
  r^{-1}\lambda_\frac13(r)$.  For $t \in[0,\tfrac13]$, let
  \begin{align*}
    \lambda_t(r) = (1-3t)r + 3t\lambda_\frac13(r)
  \end{align*}
  and for $(t,r) \in [\frac13,1]\times [0,\infty) - \{(1,0)\}$ let
  \begin{align*}
    \lambda_t(r) =
    \begin{cases}
      \lambda_\frac13(\frac{2r}{3(1-t)})& \text{$r \leq 1.5(1-t)$}\\
      2 & \text{$1.5(1-t) \leq r\leq 1.9$}\\
      \lambda_\frac13(r) & \text{$r\geq 1.9$.}
    \end{cases}
  \end{align*}
  Let $g_t(r) = (\lambda_t(r))^{-1} r$ and $g_t(0) = 1$ for $t \leq 0$
  and $g_t(0) = (1-t)$ for $0\leq t \leq 1$.  The graph of $g_t$ is
  shown in figure~\ref{fig:3} for various values of $t \in [0,1]$.
  \begin{figure}
    \centering
    \includegraphics{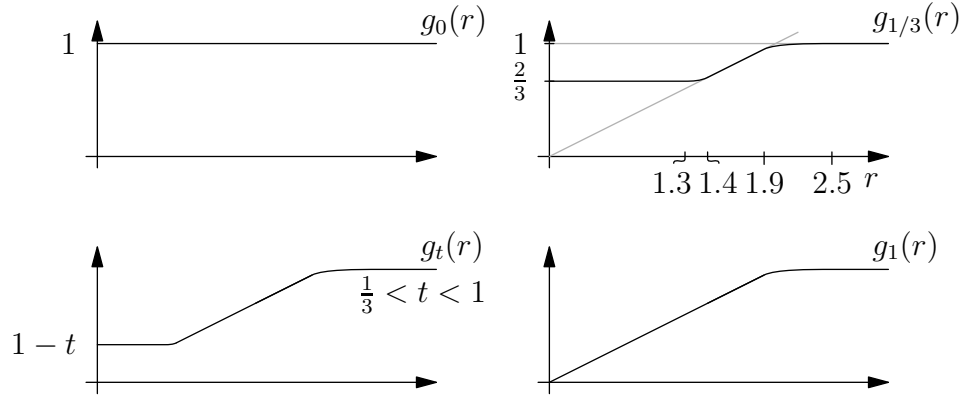}
    \caption{$g_t$ for various $t \in [0,1]$.}
    \label{fig:3}
  \end{figure}
  Define $\phi_t: \R^N \to \R^N$ by
  \begin{align*}
    \phi_t(x) = \frac{x}{g_t(|x|)}.
  \end{align*}
  Thus $\phi_t$ multiplies by $(t-1)^{-1}$ near $\phi_t^{-1}(B(0,1))$
  and is the identity outside $\phi_t^{-1}(B(0,2.5))$.  $\phi_t$
  preserves lines through the origin and $|\phi_t(x)| =
  \lambda_t(|x|)$, so the critical values of $\phi_t$ when $t \geq
  1/3$ are precisely the points in $2S^{N-1}$.  We leave $\phi_1(0)$
  undefined.
  \begin{enumerate}[(i)]
  \item For $T\subseteq G \in \Phi(\R^N)$ in collapsible position,
    define a path of subsets $G_t \subseteq \R^N$ by
    \begin{align*}
      G_t =
      \begin{cases}
        \phi_t^{-1}(G) & \text{for $t < 1$,}\\
        \{0\} \cup \phi_1^{-1}(G)& \text{for $t=1$.}
      \end{cases}
    \end{align*}
  \item For $x \in \phi_t^{-1}(T)$ or $|\phi_t(x)|\geq 3$, let $l_t(x)
    = l(\phi_t(x))$.
  \item If $x \in G_t$ has $\phi_t(x) \in B(0,3) - T$, define $l_t(x)
    = (d_t(x)-1)^2$, where $d_t$ is defined as
    \begin{align}\label{eq:49}
      d_t(x) = g_t(|x|)d(\phi_t(x))
    \end{align}
    and $d: G \cap B(0,3)\to [0,2)$ is the function from
    definition~\ref{definition:collapsible-pos}.
  \end{enumerate}
\end{construction}

The collapse of a tree $T\subseteq G$ in collapsible position in
construction~\ref{constr} above is illustrated in figure~\ref{fig:2}.
The outer gray circle in each picture is $\partial B(0,2)$, and the
region between the two gray circles is $\phi_t^{-1}(\partial
B(0,2))$.

\begin{lemma}
  For $(G,T)$ in collapsible position, the above construction gives
  elements $\Upsilon_t(G,l) = (G_t,l_t) \in \Phi(\R^N)$ for all $t \in
  [0,1]$.  Moreover, $(t,(G,l)) \mapsto (G_t,l_t)$ defines a
  continuous map $\Upsilon: [0,1] \times C \to \Phi(\R^N)$, where
  $C\subseteq \Phi(\R^N)$ is the open subspace consisting of graphs in
  collapsible position.
\end{lemma}
\begin{proof}
  The assumptions on $G$ imply that the set $\{g \in G\mid 1 \leq |g|
  \leq 3\}$ contains no vertices of $G$ and no critical points of the
  function $g \mapsto |g|$.  Therefore $\phi_t: \R^N \to \R^N$ is
  transverse to the edges of $G$ and for each vertex $v$ of $g$,
  $\phi_t$ is a diffeomorphism near $\phi_t^{-1}(v)$.  This implies
  that $G_t$ satisfies the requirements of
  definition~\ref{definition:Phi-U-as-set}, except possibly that
  parametrizations $\gamma$ satisfy $l\circ \gamma(t) = t^2$.

  Let $\gamma_t: (a,b) \to G_t$ be a parametrization of an edge of
  $G_t$ such that $\phi_t\circ \gamma_t(s)$ maps to the image of an
  incident edge, and $|\gamma_t(s)|$ is an increasing function of $s$.
  Then a direct calculation shows that $d_t(\gamma_t(s))$ has strictly
  positive derivative with respect to $s$.  Indeed, in
  formula~(\ref{eq:49}), both factors $g_t(|\gamma_t(s)|)$ and
  $d_t(\phi_t(x))$ have non-negative derivative, and at least one of
  them has strictly positive derivative.  On the subset $\phi_t^{-1}(G
  \cap \partial B(0,2))\subseteq G_t$ (the region between the gray
  circles in figure~\ref{fig:2}) the factor $d(\phi_t(\gamma_t(s)))$
  is constant and the factor $g_t(|\gamma_t(s)|)$ is $2|\gamma_t(s)|$.
  After reparametrizing $\gamma_t$ we can assume $d_t(\gamma_t(s)) = s
  +1$ in which case $\gamma_t$ is an admissible parametrization of
  $G_t$.  This implies that $(G_t,l_t) \in \Phi(\R^N)$ for all $t \in
  [0,1]$.

  Thus to each incident edge $\gamma: [0,\tau] \to G = G_{0}$
  corresponds an admissible parametrization $\gamma_t: [0,\tau] \to
  \phi_t^{-1}(\IM(\gamma))$ and $\gamma_{t'}$ and $\gamma_{t}$ have
  images that are diffeomorphic via the map $x \mapsto
  \gamma_{t'}\circ d_{t}(x)$.  These assemble to a map $G_{t} \to
  G_{t'}$ which is an isomorphism of graphs for $t,t' < 1$.  For $t' =
  1$ they assemble to a graph epimorphism $G_t \to G_1$ given by $x
  \mapsto \gamma_1\circ d_t(x)$ outside $\phi_t^{-1}(T)$ and
  collapsing $\phi_t^{-1}(T) \subseteq G_t$ to $0 \in G_1$.

  To prove continuity of $(t,(G,l)) \mapsto (G_t, l_t)$, let $t \in
  [0,1]$, $G \in C$, $u \in \R^N$.  We prove continuity at each point
  $(t,(G,l),u) \in [0,1] \times C \times \R^N$ (cf.\
  proposition~\ref{prop:sheaf-of-spaces} and the remark following its
  proof).  For $(t,u) \neq (1,0)$, continuity follows from the
  implicit function theorem, and it remains to prove continuity at
  $(1,G) \in [0,1]\times C$ at $0 \in R^N$.  This follows from the
  above mentioned graph epimorphism $G_t \to G_1$, because
  $\phi_t^{-1}(T) \subseteq B(0,1-t)$.
\end{proof}

Notice also that admissible parametrizations of $\phi_t^{-1}(G
\cap \partial B(0,2))\subseteq G_t$ will be parametrized at constant
speed.  Indeed,
\begin{align*}
  s = d_t(|\gamma_t(s)|) = g_t(|\gamma_t(s)|) d(\phi_t(\gamma_t(s))) =
  2a|\gamma_t(s)|
\end{align*}
for some constant $a = d(\phi_t(\gamma_t(s)))$.  In particular $G_1
\cap B(0,2)$ will consist of straight lines, parametrized in a linear
fashion.

If $G \in \Phi(U)$ and $e: \R^N \to U$ is an embedding such that
$e^*(G)$ is in collapsible position, we can define a path in
$\Phi(e(\R^N))$ by $t \mapsto (e^{-1})^* \circ \Upsilon_t \circ e^* (G)$.
This path is constant on $\Phi(e(\R^N - B(0,3)))$ so by the sheaf
property (proposition~\ref{prop:sheaf-of-spaces}) it glues with the
constant path $t \mapsto G|(U - e(B(0,3)))$ to a path $t \mapsto
\Upsilon_t^e(G) \in \Phi(U)$.  This defines a continuous function
$\Upsilon^e: [0,1] \times C(e) \to \Phi(U)$, where $C(e)\subseteq
\Phi(U)$ is the open subset consisting of $G$ for which $e^*(G)$ is in
collapsible position.

\begin{lemma}\label{lemma:good-disks}
  For any $(G,f) \in \Phi(M;R)$ and any $v \in \VV(R) \cap \Int M$,
  let $T_v = f^{-1}(v)$.  There exists an embedding $e = e_v: \R^N \to
  \Int M$ such that $(e^{-1}(G), e^{-1}(T_v))$ is collapsible.  If
  $W\subseteq U$ is a neighborhood of $T_v$ then $e_v$ can be chosen
  to have image in $W$.  In particular we can choose the embeddings
  $e_v, v \in \VV(R)\cap M$ to have disjoint images.
\end{lemma}
\begin{proof}
  Embed small disks around each vertex of $T_v$, and do connected sum
  along a small tubular neighborhood of each edge of $T_v$.
\end{proof}

If the embeddings in the above lemma have disjoint images, we get an
embedding $e: (\VV(R) \cap \Int M) \times \R^N \to \Int M$.  We will
say that $e$ and $(G,f)$ are \emph{compatible} if they satisfy the
conclusion of the lemma: $(e_v^{-1}(G), e_v^{-1}(T_v))$ is collapsible
for all $v \in (\VV(R)\cap \Int M)$, where $e_v = e(v,-): \R^N \to
\Int M$.  Thus the lemma says that for any $(G,f) \in \Phi(M;R)$ we
can find arbitrarily small compatible embeddings $e$.

From a compatible embedding $e: (\VV(R) \cap \Int M) \times \R^N \to \Int
M$ we construct a path $t \mapsto \Upsilon_t^e(G) \in \Phi(U)$ as
above, i.e.\ by gluing the path
\begin{align*}
  t \mapsto \prod_v (e_v^{-1})^* \circ \Upsilon_t \circ e_v^*(G) \in
  \Phi(\coprod_v e_v(\R^N))
\end{align*}
with the constant path
\begin{align*}
  t \mapsto G |(U - \coprod_v e_v(\R^N)).
\end{align*}
If $e$ is compatible with $(G,f)$, the path $t \mapsto \Upsilon^e_t(G)
\in \Phi(U)$ has a unique lift to a path $[0,1] \to \Phi(M;R)$ which
starts at $(G,f)$.  We will use the same notation $t \mapsto
\Upsilon_t^e(G,f)$ for the lifted path $[0,1]\to \Phi(M;R)$.  We point
out that $\Upsilon^e_t(G,f) \in \Emb_R(M)$ if $(G,f) \in \Emb_R(M)$
and that $\Upsilon^e_1(G,f) \in \Emb_R(M)$ for all $(G,f) \in
\Phi(M;R)$.  If we let $C(e) \subseteq \Phi(M;R)$ be the set of
$(G,f)$ which are compatible with $e$, we have constructed a
\emph{homotopy} between the inclusion $C(e) \subseteq \Phi(M;R)$ and
the map $\Upsilon_1^e: C(e) \to \Emb_R(M)$.

We can now prove that $\pi_0(\Phi(M;R),\Emb_R(M)) = 0$.  Namely, let
$(G,f) \in \Phi(M;R)$ and use lemma~\ref{lemma:good-disks} to choose a
compatible $e: (\VV(R) \cap \Int M) \times \R^N \to \Int M$.  Then
$t \mapsto \Upsilon_t^e(G,f)$ is a path to a point in $\Emb_R(M)$ as
required.

To prove that the higher relative homotopy groups vanish, we need to
carry out the above collapsing in families, i.e.\ when parametrized by
a map $X \to \Phi(M;R)$.  Unfortunately it does not seem easy to prove
a pa\-ra\-me\-trized version of lemma~\ref{lemma:good-disks}.  Instead
we use collapsing along multiple $e: \R^N \to U$ at once.
\begin{definition}
  Let $Q$ be the set of all embeddings $e: (\VV(R) \cap M) \times \R^N
  \to M$.  Write $e < e'$ if $e(\{v\} \times \R^N)\subseteq e'(\{v\}
  \times B(0,1))$ for all $v \in \VV(R)$.  This makes $Q$ into a
  poset.  We give $Q$ the discrete topology.  Let $P\subseteq Q \times
  \Phi(M;R)$ be the subspace consisting of $(e,(G,f))$ with $(G,f) \in
  C(e)$.  $P$ is topologized in the product topology and ordered in
  the product ordering, where $\Phi(M;R)$ has the trivial order.
\end{definition}
\begin{lemma}\label{lemma:map-of-pairs-is-eq}
  The projection $p: BP \to \Phi(M;R)$ is a weak equivalence.  The
  restriction to $p^{-1}(\Emb_R(M)) \to \Emb_R(M)$ is also a weak
  equivalence.
\end{lemma}
The proof is based on lemma~\ref{lemma:Segal}.  First recall that any
poset $D$ can be considered as a category.  It is easy to see that
$BD$ is contractible when it has a subset $C\subseteq D$ which in the
induced ordering is totally ordered and \emph{cofinal}, i.e.\ if for
every $d,d' \in D$ there is $c \in C$ with $d \leq c$ and $d' \leq c$.
Indeed, any finite subcomplex of $\|N_\bullet D\|$ will be contained
in the star of some vertex $c \in \|N_\bullet C\| \subseteq
\|N_\bullet D\|$.
\begin{proof}[Proof of lemma~\ref{lemma:map-of-pairs-is-eq}]
  $p$ is induced by the projection $\pi: P \to \Phi(M;R)$ which is
  etale (because $N_kP \subseteq N_kQ \times \Phi(M;R)$ is open and
  $N_kQ$ is discrete), so by lemma~\ref{lemma:Segal} it suffices to
  prove that $B(\pi^{-1}(G,f))$ is contractible for all $(G,f)$.  For
  any $(G,f)$ we can choose, by lemma~\ref{lemma:good-disks}, a
  sequence $e_n: (\VV(R) \cap \Int M) \times \R^N \to \Int M$, $n \in
  \N$ of embeddings compatible with $(G,f)$, such that $e_1 > e_2 >
  \dots$, and with $e_n(\{v\} \times \R^N)$ contained in the
  $(1/n)$-neighborhood of $T_v$.  This totally ordered subset of
  $\pi^{-1}((G,f))$ is cofinal.  The second part is proved the same
  way.
\end{proof}

\begin{definition}
  For $t = (t_0, \dots, t_k) \in [0,1]^{k+1}$ and $\chi = (e_0 < \dots
  < e_k,(G,f)) \in N_k(P)$, let
  \begin{align*}
    \Upsilon((t_0, \dots, t_k), \chi) = \Upsilon_{t_k}^{e_k} \circ
    \dots \circ \Upsilon_{t_0}^{e_0} (G,f).
  \end{align*}
  This defines a continuous map $\Upsilon: [0,1]^{k+1} \times N_k(P)
  \to \Phi(M;R)$.
\end{definition}

\begin{proposition}\label{prop:collaps-gives-we}
  The inclusion $\Emb_R(M) \to \Phi(M;R)$ is a weak equivalence.
\end{proposition}
\begin{proof}
  Let $m: \Delta^k \to [0,1]^{k+1}$ be defined by
  \begin{align*}
    m(t_0, \dots, t_k) = (t_0, \dots, t_k)/\max(t_0, \dots, t_k).
  \end{align*}
  Then the maps $h: [0,1] \times \Delta^k \times N_k(P) \to \Phi(M;R)$
  defined by
  \begin{align*}
    h(\tau,t,\chi) = \Upsilon(\tau m(t), \chi)
  \end{align*}
  glue together to a map $h:[0,1] \times BP \to \Phi(M;R)$ which is a
  homotopy between the projection map $p: BP \to \Phi(M;R)$ and the
  map
  \begin{align*}
    q = h(1,-): BP \to \Emb_R(M) \subseteq \Phi(M;R).
  \end{align*}
  This produces a null homotopy of the map of pairs
  \begin{align*}
    p: (BP, p^{-1}(\Emb_R(M))) \to (\Phi(M;R), \Emb_R(M)),
  \end{align*}
  which together with lemma~\ref{lemma:map-of-pairs-is-eq} proves that
  $\pi_*(\Phi(M;R), \Emb_R(M)) = 0$.
\end{proof}

\begin{proposition}\label{prop:emb-connectivity}
  The space $\Emb_R(M)$ is $(N-4)$-connected when $\Int M$ is
  $(N-3)$-connected.
\end{proposition}
\begin{proof}[Proof sketch.]
  Thinking about $\Emb_R(M)$ as a space of embeddings $j: R \to U$
  this is mostly standard, although the presence of vertices deserves some
  comment.

  Firstly, we can fix $j$ on $\VV(R) \cap \Int M$, since this changes
  homotopy groups only in degrees above $(N-4)$.  Secondly, the proof
  of proposition~\ref{prop:collaps-gives-we} shows that we can assume
  there is a ball $B = B(v,\epsilon_v)$ around each vertex $v \in
  \VV(R) \cap \Int M$ such that $j$ is linear on $j^{-1}(B)$ (the
  point is that by construction, this holds after applying
  $\Upsilon^e_1$).  Thirdly, we can fix $j$ on $j^{-1}B$, since this
  changes homotopy groups only in degrees above $(N-3)$.  Then we are
  reduced to considering embeddings of a disjoint union of intervals
  into the manifold $\Int M - \cup_v \Int B(v,\epsilon_v)$, and these
  form an $(N-4)$-connected space when $\Int M$ is $(N-3)$-connected
  (an easy consequence of transversality).
\end{proof}

\begin{proof}[Proof of theorem~\ref{thm:Phi-is-B-Cat}]
  We have the maps
  \begin{align*}
    B\mathcal{G}_S \leftarrow \hocolim_{R \in \mathcal{G}_S} \Phi(M;R)
    \rightarrow \Phi^S(M).
  \end{align*}
  The map pointing to the right is a weak equivalence by
  proposition~\ref{prop:hocolim-weak-equiv}.
  The map pointing to the left is obtained by taking $\hocolim$ of the
  collapse map
  \begin{align*}
    \Phi(M;R) \to \mathit{point}
  \end{align*}
  which is $(N-3)$-connected by
  propositions~\ref{prop:collaps-gives-we} and
  \ref{prop:emb-connectivity}.
\end{proof}

\subsection{Abstract graphs}
\label{sec:abstract-graphs}

The goal in this section is to determine the homotopy type of the
space $B\mathcal{G}_S$.  Although the objects of $\mathcal{G}_S$ are
embedded graphs, the embeddings play no role in the morphisms, and the
category $\mathcal{G}_S$ is equivalent to a combinatorially defined
category of \emph{abstract} graphs.  We recall the definition of
abstract graphs, cf.\ \cite{MR721773}.

\begin{definition}
  \begin{enumerate}[(i)]
  \item A finite \emph{abstract graph} is a finite set $G$ with an
    involution $\sigma: G \to G$ and a retraction $t: G \to G^\sigma$
    onto the fixed point set of $\sigma$.
  \item The \emph{vertices} of $G$ is the set $G^\sigma$ of fixed
    points of $\sigma$, and the complement $G - G^\sigma$ is the set
    of \emph{half-edges}.  The \emph{valence} of a vertex $x \in
    G^\sigma$ is the number $v(x) = |t^{-1}(x)| - 1$.  In this paper,
    all graphs are assumed not to have vertices of valence $0$ and
    $2$.
  \item A \emph{leaf} of $G$ is a valence 1 vertex.  A \emph{leaf
      labelling} of $G$ is an identification of the set of leaves of
    $G$ with $\{1,\dots, s\}$.
  \item A \emph{cellular map} $G \to G'$ between two abstract graphs
    is a set map preserving $\sigma$ and $t$.
  \item A cellular map $G \to G'$ is a \emph{graph epimorphism} if the
    inverse image of each half-edge of $G'$ is a single half-edge of
    $G$, and the inverse image of a vertex of $G'$ is a \emph{tree}
    (i.e.\ contractible graph), not containing any leaves.  If the
    leaves of $G$ and $G'$ are labelled, we require the map to
    preserve the labelling.
  \item For $s \geq 0$, let $\mathcal{G}_s$ denote the category whose
    objects are finite abstract graphs with leaves labelled by
    $\{1,\dots, s\}$ and whose morphisms are graph epimorphisms.
  \end{enumerate}
\end{definition}

A cellular map is a graph epimorphism if and only if it can be written
as a composition of \emph{isomorphisms} and \emph{elementary
  collapses}, i.e.\ maps which collapse a single non-loop, non-leaf
edge to a point.

\begin{lemma}\label{lemma:abstract-graphs}
  Let $N \geq 3$ and let $M\subseteq U$ be compact with $\Int M$
  connected.  Let $S \in \Phi(U - \Int M)$ be a germ of a graph with
  $s\geq 0$ ends in $\Int M$ (i.e.\ $s$ is the cardinality of the
  inverse limit of $\pi_0(S \cap (\Int M-K))$ over larger and larger
  compact sets $K \subseteq \Int M$).  Then we have an equivalence of
  categories
  \begin{align}
    \mathcal{G}_S \simeq \mathcal{G}_s.\tag*{\qed}
  \end{align}
\end{lemma}

It remains to determine the homotopy type of the space $B
\mathcal{G}_s$.  A finite abstract graph $G$ has a \emph{realization}
\begin{align*}
  |G| = \big( G^\sigma \amalg ((G - G^\sigma) \times [-1,1]) \big) /
  \sim,
\end{align*}
where $\sim$ is the equivalence relation generated by $(x,r) \sim
(\sigma x, -r)$ and $(x,1) \sim t(x) \in G^\sigma$ for $x \in G -
G^\sigma$, $r \in [-1,1]$.  Let $\partial |G| \subseteq |G|$ be the
set of leaves (valence 1 vertices).  Let $\Aut(G)$ denote the group of
homotopy classes of homotopy equivalences $|G| \to |G|$ restricting to
the identity on $\partial |G|$.  Recall from
section~\ref{subsec:outline-proof} that $A_n^s = \Aut(G_n^s)$, where
$G_n^s$ is a graph with first Betti number $n$ and $s$ leaves.  In
particular $A_n^0 = \Out(F_n)$ and $A_n^1 = \Aut(F_n)$.

$B\mathcal{G}_0$ and $B\mathcal{G}_1$ are directly related to
automorphisms of free groups, via Culler-Vogtmann's \emph{outer space}
\cite{MR830040}.  As mentioned in section~\ref{subsec:outline-proof},
outer space is a contractible space with an action of $\Out(F_n)$.
Culler-Vogtmann also define a certain subspace called the \emph{spine
  of outer space}, which is an equivariant deformation retract.  For a
finite abstract graph $G_0$, the spine of outer space $X(G_0)$ has one
simplex for each isomorphism class of pairs $(G,h)$, where $G$ is a
finite abstract graph and $h: |G| \to |G_0|$ is a homotopy class of a
homotopy equivalence.  $(G,h)$ is a face of $(G',h')$ when there
exists a graph epimorphism $\phi: G' \to G$ such that $h \circ |\phi|
\simeq h'$.  Culler-Vogtmann prove that $X(G_0)$ is contractible.
(They state this only for connected $G_0$; the general statement
follows from the homeomorphism $X(G_0 \amalg G_1) \cong X(G_0) \times
X(G_1)$.)
\begin{proposition}\label{prop:igusa}
  There is a homotopy equivalence
  \begin{align}\label{eq:20}
    B\mathcal{G}_s \simeq \coprod_{G} B \Aut(G),
  \end{align}
  where the disjoint union is over finite graphs $G$ with $s$ leaves,
  one of each homotopy type.
\end{proposition}
The right hand side of the homotopy equivalence~\eqref{eq:20} can
conveniently be reformulated in terms of a category
$\mathcal{G}_s^\simeq$.  The objects of $\mathcal{G}_s^\simeq$ are the
objects of $\mathcal{G}_s$, but morphisms $G \to G'$ in
$\mathcal{G}_s^\simeq$ are homotopy classes of homotopy equivalences
$(|G|, \partial |G|) \to (|G'|, \partial |G'|)$, compatible with the
labellings.  We have inclusion functors
\begin{align*}
  \xymatrix{
    {\mathcal{G}_s} \ar[r]^-{f} & {\mathcal{G}_s^\simeq}&
    {\coprod_{G} \Aut(G)} \ar[l]_-{g}
  }
\end{align*}
where $f$ is the identity on the set objects and takes geometric
realization of morphisms, and $g$ is the inclusion of a skeletal
subcategory.  Consequently $g$ is an equivalence of categories, and
the statement of proposition~\ref{prop:igusa} is equivalent to $f$
inducing a homotopy equivalence $Bf: B\mathcal{G}_s \to
B\mathcal{G}_s^\simeq$.

\begin{proof}[Proof sketch]
  We first consider the case $s = 0$, following \cite[Theorem
  8.1.21]{MR1945530}.

  For a fixed object $G_0 \in \mathcal{G}^\simeq$, we consider the
  over category $(\mathcal{G}\downarrow G_0)$.  Its objects are pairs
  $(G,h)$ consisting of an object $G \in \Ob(\mathcal{G})$ and a
  homotopy class of a homotopy equivalence $h: |G| \to |G_0|$.  Its
  morphisms $(G,h) \to (G', h')$ are graph epimorphisms $\phi: G \to
  G'$ with $h' \circ |\phi| \simeq h$.  It is equivalent to the
  opposite of the poset of simplices in the spine of outer space, and
  hence contractible.  Then the claim follows from Quillen's ``theorem
  A'' (\cite{MR0338129}).

  We proceed by induction in $s$.  Recall that $\Aut(G) =
  \pi_0h\Aut(G)$, where $h\Aut(G)$ is the topological monoid of
  self-homotopy equivalences of $|G|$ restricting to the identity on
  the boundary.  Every connected component of $h\Aut(G)$ is
  contractible and we have $Bh\Aut(G) \simeq B\Aut(G)$.  The monoid
  $h\Aut(G)$ acts on $|G|$, and the Borel construction is
  \begin{align*}
    Eh\Aut(G) \times_{h\Aut(G)} |G| \simeq \coprod_{p} B h\Aut(G'),
  \end{align*}
  where $G'$ is obtained by attaching an extra leaf to $G$ at a point
  $p$.  The disjoint union is over $p \in |G| - \partial |G|$, one in
  each $h\Aut(G)$-orbit.  It follows that the map
  \begin{align*}
    B\mathcal{G}_{s+1}^\simeq \to B\mathcal{G}_s^\simeq,
  \end{align*}
  induced by forgetting the leaf labelled $s+1$, has homotopy fiber
  $|G|$ over the point $G \in B\mathcal{G}_s^\simeq$.

  Let $\Gamma: \mathcal{G}_s \to \mathrm{CAT}$ be the functor which to
  $G\in \mathcal{G}_s$ associates the poset of simplices of $G$ which
  are not valence 1 vertices, ordered by reverse inclusion.  Recall
  from e.g.\ \cite{MR510404} that to such a functor there is a
  associated category $(\mathcal{G}_s \wr \Gamma)$.  An object of the
  category $(\mathcal{G}_s \wr \Gamma)$ is a pair $(G, \sigma)$, with
  $G \in \mathcal{G}_s$ and $\sigma \in \Gamma(G)$ and a morphism
  $(G,\sigma) \to (G', \sigma')$ is a pair $(\phi,\psi)$ with $\phi: G
  \to G'$ and $\psi: \Gamma(\phi)(\sigma) \to \sigma'$.  There is a
  functor
  \begin{align*}
    (\mathcal{G}_s \wr \Gamma) \to \mathcal{G}_{n+1}
  \end{align*}
  which maps $(G,\sigma)$ to the graph obtained by attaching a leaf
  labeled $s+1$ to $G$ at the barycenter of $\sigma$.  This is an
  equivalence of categories, and it follows (by \cite{MR510404}) that
  the homotopy fiber of the projection $B \mathcal{G}_{s+1} \to B
  \mathcal{G}_s$ over the point $G \in B \mathcal{G}_s$ is
  $B(\Gamma(G)) \cong |G|$.

  Therefore the diagram
  \begin{align*}
    \xymatrix{
      {B\mathcal{G}_{s+1}} \ar[r]\ar[d] & {B\mathcal{G}_{s+1}^\simeq}\ar[d] 
      \\
      {B\mathcal{G}_s} \ar[r] & {B\mathcal{G}_s^\simeq}
    }
  \end{align*}
  is homotopy cartesian.  This proves the induction step.
\end{proof}

Summarizing theorem~\ref{thm:Phi-is-B-Cat},
lemma~\ref{lemma:abstract-graphs}, and proposition~\ref{prop:igusa} we
get
\begin{theorem}
  \label{thm:conclusion-htpy-types}
  Let $N \geq 3$, let $U\subseteq \R^N$ be open, and let $M\subseteq
  U$ be compact with $\Int M$ $(N-3)$-connected.  Let $S \in \Phi(U -
  \Int M)$ be a germ of a graph with $s$ ends in $\Int M$.  Then we
  have an $(N-3)$-connected map
  \begin{align}
    \Phi^S(M) \to \coprod_{G} B \Aut(G),
  \end{align}
  where the disjoint union is over finite graphs $G$ with $s$ leaves,
  one of each homotopy type.
\end{theorem}

\subsection{$B\Out(F_n)$ and the graph spectrum}
\label{subsec:param-pontrj-thom}

We are now ready to begin the proof outlined in
subsection~\ref{subsec:outline-proof}.  The first goal is to define
the maps~\eqref{eq:2} and \eqref{eq:3}.  The space $B_N$ in the
following definition is the domain of the map~\eqref{eq:2}.

\begin{definition}
  \label{defn:B-N}
  Let $I = [-1,1]$.  Let $B_N \subseteq \Phi(\R^N)$ be the subset
  \begin{align*}
    B_N = \Phi^{[\emptyset]}(I^N),
  \end{align*}
  i.e.\ the set of graphs contained in $\Int(I^N)$.
\end{definition}

The homotopy type of the space $B_N$ is determined by
theorem~\ref{thm:conclusion-htpy-types}.

\begin{proposition}\label{prop:B-N-htpy-type}
  There is an $(N-3)$-connected map
  \begin{align*}
    B_N \to \coprod_G B \Aut(G),
  \end{align*}
  where the disjoint union is over graphs $G$ without leaves, one of
  each homotopy type.  Consequently we have a weak equivalence
  \begin{multline*}
    \{G \in B_\infty \mid \text{there \emph{exists} a homotopy
      equivalence $G \simeq \vee^n S^1$} \} \simeq\\
    B\Out(F_n).\tag*{\qed}
  \end{multline*}
\end{proposition}

Approximating $B\Out(F_n)$ by the space consisting of $G \in B_N$ for
which there exists a homotopy equivalence $G \simeq \vee^n S^1$ is
analogous to the approximation
\begin{align*}
  B\Diff(M) \sim \Emb(M, \Int (I^N))/ \Diff(M)
\end{align*}
for a smooth manifold $M$.  The right hand side is the space of
submanifolds $Q \subseteq \Int (I^N)$ for which there exists a
diffeomorphism $Q \cong M$.

The empty set $\emptyset \subseteq \R^N$ is a graph, and we consider
it the basepoint of $\Phi(\R^N)$.

\begin{definition}
  Let $\epsilon_N: S^1 \wedge \Phi(\R^N) \to \Phi(\R^{N+1})$ be the
  map induced by the map $\R \times \Phi(\R^N) \to \Phi(\R^{N+1})$
  given by $(t,G) \mapsto \{-t\} \times G$.
\end{definition}
\begin{lemma}
  $\epsilon_N$ is well defined and continuous.
\end{lemma}
\begin{proof}
  For any compact subset $K \subseteq \R^N$ with $K \subseteq
  cD^{N+1}$ we will have $(\{t\} \times G)\cap K = \emptyset \in
  \Phi(\R^{N+1})$ as long as $|t| > c$ or $G \cap cD^N = \emptyset$.
  This proves that $\epsilon_N$ is continuous at the basepoint.
  Continuity on $\R \times \Phi(\R^N)$ follows from
  proposition~\ref{prop:res-is-cont}.
\end{proof}

\begin{definition}
  Let $\mathbf{\Phi}$ be the spectrum with $N$th space $\Phi(\R^N)$
  and structure maps $\epsilon_N$.  This is the \emph{graph spectrum}.
\end{definition}
We will not use any theory about spectra.  In fact we will always work
with the corresponding infinite loop space $\Omega^\infty
\mathbf{\Phi}$ defined as
\begin{align*}
  \Omega^\infty\mathbf{\Phi} = \colim_{N \to \infty} \Omega^N
  \Phi(\R^N)
\end{align*}
where the map $\Omega^N \Phi(\R^N) \to \Omega^{N+1} \Phi(\R^{N+1})$ is
the $N$-fold loop of the adjoint of $\epsilon_N$.

$\mathbf{\Phi}$ is the analogue for graphs of the spectrum $\MT{d}$
for $d$-manifolds in the paper \cite{math.AT/0605249}.  The analogy is
clarified in chapter~\ref{sec:some-remarks-manif}, especially
proposition~\ref{prop:6.2}. $\MT{d}$ is the Thom spectrum of the
universal stable normal bundle for $d$-manifold bundles, $-U_d \to
BO(d)$.  Thus $\mathbf{\Phi}$ is a kind of ``Thom spectrum of the
universal stable normal bundle for graph bundles.''
Remark~\ref{remark:generalized-spherical-fib} explains in what sense
$\mathbf{\Phi}$ is the Thom spectrum of a ``generalized stable
spherical fibration''.  In this subsection we will define a map which,
alluding to a similar analogy, we could call the \emph{parametrized
  Pontryagin-Thom collapse map for graphs}
\begin{align}
  \label{eq:5}
  B\Out(F_n) \to \Omega^\infty\mathbf{\Phi}.
\end{align}

Given $G \in B_N$ and $v \in \R^N$ we can translate $G$ by $v$ and get
an element
\begin{align*}
  \tau_N(G)(v) = G - v \in \Phi(\R^N).
\end{align*}
We have $\tau_N(G)(v) \to \emptyset$ if $|v| \to \infty$, so $\tau_N$
extends uniquely to a continuous map
\begin{align}\label{eq:37}
  ({B_N}) \wedge S^N \xrightarrow{\tau_N} \Phi(\R^N).
\end{align}
\begin{definition}
  \label{defn:alpha}
  Let $\tau_N: B_N \to \Omega^N \Phi(\R^N)$ be the adjoint of the
  map~\eqref{eq:37}.
\end{definition}

In the following diagram, the left vertical map $B_N \to B_{N+1}$ is
the inclusion $G \mapsto \{0\} \times G$.
\begin{align*}
  \xymatrix{ B_N \ar[r]^-{\tau_N} \ar[d] & {\Omega^N \Phi(\R^N)}
    \ar[d]^{\epsilon_N} \\
    B_{N+1} \ar[r]^-{\tau_{N+1}} & {\Omega^{N+1} \Phi(\R^{N+1})}.}
\end{align*}
The diagram is commutative, and we get an induced map
\begin{align}\label{eq:38}
  \tau_\infty : B_\infty \to \Omega^\infty\mathbf{\Phi}.
\end{align}
By theorem~\ref{thmcor:main-decomposition}, $B\Out(F_n)$ is a
connected component of $B_\infty$, and we define the map~\eqref{eq:5}
as the restriction of~\eqref{eq:38}.

The map $\tau_N$ is homotopic to a map $\tilde \tau_N$ defined in a
different way.  This construction will be used in
section~\ref{subsubsec:second-proof}, but is not logically necessary
for the proof of theorem~\ref{thm:main}.  $\tilde \tau_N$ is similar
to the ``scanning'' map of
\cite{MR533892}, and is defined as follows.  Choose a map
\begin{align*}
  e: \Int(I^N) \times \R^N = T\Int(I^N) \to\Int(I^N)
\end{align*}
such that for each $p \in \Int(I^N)$, the induced map
\begin{align*}
  e_p: T_p\Int(I^N) \to \Int(I^N).
\end{align*}
is an embedding with $e_p(0) = p$ and $De_p(0) = \id$.  We can arrange
that the radius of $e_p(T_p\Int(I^N))$ is smaller than the distance
$\dist(p,\partial I^N)$.  To a graph $G \in B_N$ and a point $p \in
\Int(I^N)$ we associate
\begin{align*}
  \tilde \tau(G)(p) = (e_p)^*(G) \in \Phi(T_p\R^N) = \Phi(\R^N).
\end{align*}
By proposition~\ref{prop:res-is-cont}, the action of $\Diff(\R^k)$ on
$\Phi(\R^k)$ is continuous.  Therefore we can apply $\Phi$ fiberwise
to vector bundles: If $V \to X$ is a vector bundle, then there is a
fiber bundle $\Phi^\mathrm{fib}(V)$ whose fiber over $x$ is
$\Phi(V_x)$.  We will have $\tilde \tau_N(G)(p) = \emptyset$ for all
$p$ outside some compact subset of $\Int(I^N)$.  $\tilde \tau_N(G)(p)$
is continuous as a function of $p$, and can be interpreted as a
\emph{compactly supported section} over $\Int(I^N)$ of the fiber
bundle $\Phi^\mathrm{fib}(T \R^N)$.  We define $\tilde \tau_N(G)(p) =
\emptyset$ for $p \in \partial I^N$ and get a section
\begin{align*}
  \tilde \tau_N(G) \in \Gamma((I^N,\partial
  I^N),\Phi^\mathrm{fib}(T\R^N)) \cong \Omega^N\Phi(\R^N).
\end{align*}
$\tilde \tau_N(G)$ depends continuously on $G$ so we get a continuous
map
\begin{align}\label{eq:39}
  \tilde \tau_N: B_N \to \Omega^N \Phi(\R^N),
\end{align}
which is easily seen to be homotopic to the map $\tau_N$ of
definition~\ref{defn:alpha}.

\section{The graph cobordism category}
\label{sec:graph-cobordism-category}

\begin{sloppypar}
  As explained in the introduction, composing~\eqref{eq:5} with the
  maps $B\Aut(F_n) \to B\Aut(F_{n+1}) \to B\Out(F_{n+1})$ gives a map
\begin{align}\label{eq:7}
  \coprod_{n\geq 0} B \Aut(F_n) \xrightarrow{\tau}
  \Omega^\infty\mathbf{\Phi}.
\end{align}
We will prove the following.
\begin{theorem}\label{theorem:main-of-sect-3}
  $\tau$ induces a homology equivalence
  \begin{align*}
    \Z \times B \Aut_\infty \to \Omega^\infty \mathbf{\Phi}.
  \end{align*}
\end{theorem}
\end{sloppypar}
$\coprod B\Aut(F_n)$ is a topological \emph{monoid} whose group
completion is $\Z \times B\Aut_\infty^+$.  It turns out to be fruitful
to enlarge it to a topological \emph{category}, with more than one
object.  Namely we will define a ``graph cobordism category''
$\mathcal{C}_N$ whose morphisms are graphs in $\R^N$.
\begin{definition}
  \label{defn:graph-cobordism-cat}
  For $\epsilon > 0$, let $\Ob(\mathcal{C}_N^\epsilon)$ be the set
  \begin{align*}
    \{(a,A,\lambda) \mid \text{$a \in \R$, $A \subseteq \Int(I^{N-1})$
      finite, $\lambda\in (-1+\epsilon,1-\epsilon)^A$}\}.
  \end{align*}
  For an object $c = (a,A,\lambda)$, let $U^\epsilon_a = (a-\epsilon,
  a+\epsilon)\times \R^{N-1}$ and let
  \begin{align*}
    S^\epsilon_{c} = (a-\epsilon, a+\epsilon) \times A.
  \end{align*}
  Equipped with the map $l_c: S_c \to [0,1)$ given by
  \begin{align*}
    l_c(a + t,x) = (t + \lambda(x))^2
  \end{align*}
  for $|t| < \epsilon$, this defines an element $(S_c^\epsilon, l_c)
  \in \Phi(U_a^\epsilon)$.  For two objects $c_0 =
  (a_0,A_0,\lambda_0)$ and $c_1 = (a_1,A_1,\lambda_1)$ with $0 <
  2\epsilon < a_1 - a_0$, let $\mathcal{C}_N^\epsilon(c_0, c_1)$ be
  the set consisting of $(G,l) \in \Phi((a_0 - \epsilon, a_1 +
  \epsilon) \times \R^{N-1})$ satisfying $(G,l) | U^\epsilon_{a_\nu} =
  (S^\epsilon_{c_\nu},l_{c_\nu})$ for $\nu = 0,1$.  If $c_2 = (a_2,
  A_2, \lambda_2)$ is a third object and $(G',l') \in
  \mathcal{C}^\epsilon_N(c_1, c_2)$, let $(G,l) \circ (G',l') = (G'',
  l'')$, where
  \begin{align*}
    G'' = G \cup G'
  \end{align*}
  and $l'': G'' \to [0,1]$ agrees with $l$ on $G$ and with $l'$ on
  $G'$.  This defines $\mathcal{C}_N^\epsilon$ as a category of sets.
  Topologize the total set of morphisms as a subspace of
  \begin{align*}
    \coprod_{a_0, a_1} \Phi((a_0 - \epsilon, a_1 + \epsilon) \times
    \R^{N-1}),
  \end{align*}
  where the coproduct is over $a_0, a_1\in \R$ with either $a_0 = a_1$
  (the identities) or $0 < 2\epsilon < a_1 - a_0$.  We have inclusions
  $\mathcal{C}_N^\epsilon \to \mathcal{C}_N^{\epsilon'}$ when
  $\epsilon' < \epsilon$, and we let
  \begin{align*}
    \mathcal{C}_N = \colim_{\epsilon \to 0} \mathcal{C}^\epsilon_N.
  \end{align*}
\end{definition}

The following theorem determines the homotopy type of the space of
morphisms between two fixed objects in $\mathcal{C}_N$.  It is a
consequence of theorem~\ref{thm:Phi-is-B-Cat},
lemma~\ref{lemma:abstract-graphs}, and proposition~\ref{prop:igusa}.

\begin{theorem}\label{thmcor:main-decomposition}
  Let $c_0 = (a_0, A_0, \lambda_0)$ and $c_1 = (a_1, A_1, \lambda_1)$
  be objects of $\mathcal{C}_N$ with $a_0 < a_1$.  There is an
  $(N-3)$-connected map
  \begin{align*}
    \mathcal{C}_N(c_0, c_1) \to \coprod_{G} B \Aut(G),
  \end{align*}
  where the disjoint union is over finite graphs $G$ with $s = |A_0| +
  |A_1|$ leaves, one of each homotopy type.  Consequently
  \begin{align}\label{eq:29}
    \{G \in \mathcal{C}_\infty(c_0,c_1) \mid \text{$G$ is connected}\}
    \simeq \coprod_{n \geq 0} BA_n^{s}.
  \end{align}
  ($n=0$ should be excluded if $s = 1$ and $n= 0,1$ should be excluded
  if $s = 0$.)\qed
\end{theorem}

For the proof of theorem~\ref{theorem:main-of-sect-3} we need two more
definitions.
\begin{definition}
  Let $D_N\subseteq \Phi(\R^N)$ denote the subspace
  \begin{align*}
    \{G \in \Phi(\R^N) \mid G \subseteq \R \times \Int(I^{N-1}) \}.
  \end{align*}
\end{definition}
\begin{definition}\label{defn:pos-bd}
  The \emph{positive boundary} subcategory
  $\mathcal{C}_N^\partial\subseteq \mathcal{C}_N$ is the subcategory
  with the same space of objects, but whose space of morphisms from
  $c_0 = (a_0,A_0,\lambda_0)$ to $c_1 = (a_1,A_1,\lambda_1)$ is the
  subset
  \begin{align*}
    \{ G \in \mathcal{C}_N(c_0, c_1) \mid \text{$A_1 \to \pi_0(G)$
      surjective} \}.
  \end{align*}
\end{definition}

Then theorem~\ref{theorem:main-of-sect-3} is proved in the following
four steps.  Carrying them out occupies the remainder of this chapter.
\begin{itemize}
\item There is a homology equivalence
  \begin{align}\label{eq:53}
    \Z \times B \Aut_\infty \to \Omega B \mathcal{C}_\infty^\partial.
  \end{align}
\item There is a weak equivalence
  \begin{align}\label{eq:54}
    B \mathcal{C}_N \simeq D_N.
  \end{align}
\item There is a weak equivalence
  \begin{align}\label{eq:55}
    D_N \xrightarrow{\simeq} \Omega^{N-1}\Phi(\R^N).
  \end{align}
\item The inclusion induces a weak equivalence
  \begin{align}\label{eq:56}
    B\mathcal{C}_\infty^\partial \xrightarrow{\simeq} B \mathcal{C}_\infty.
  \end{align}
\end{itemize}

Then theorem~\ref{theorem:main-of-sect-3} follows by looping
(\ref{eq:54}), (\ref{eq:55}), and (\ref{eq:56}), taking the direct
limit $N \to \infty$ in (\ref{eq:54}) and (\ref{eq:55}), and
composing.

\subsection{Poset model of the graph cobordism category}
\label{subsec:poset-model-graph}

We will use $\R^\delta$ to denote the set $\R$ of real numbers,
equipped with the \emph{discrete} topology.

\begin{definition}
  \begin{enumerate}[(i)]
  \item Let $D_N^\pitchfork\subseteq \R^\delta \times D_N$ be the
    space of pairs $(a,(G,l))$ satisfying
    \begin{align}
      \label{eq:18}
      G \pitchfork \{a\}\times \R^{N-1}.
    \end{align}
    This is a poset, with ordering defined by $(a_0,G) \leq (a_1,G')$
    if and only if $G = G'$ and $a_0 \leq a_1$.
  \item For $\epsilon > 0$, let $D_N^{\perp,\epsilon} \subseteq
    D_N^\pitchfork$ be the subposet defined as follows.  $(a,(G,l))
    \in D_N^{\perp,\epsilon}$ if there exists $c = (a,A,\lambda)$ as
    in definition~\ref{defn:graph-cobordism-cat} such that
    $(G,l)|U^\epsilon_a = (S_c^\epsilon,l_c)$.
  \item Let $D_N^\perp$ be the colimit of $D_N^{\perp,\epsilon}$ as
    $\epsilon \to 0$.
  \end{enumerate}
\end{definition}

There is an inclusion functor $i: D_N^\perp \to D_N^\pitchfork$, and a
forgetful map $u: D_N^\pitchfork \to D_N$.  There is also a functor
$c: D_N^\perp \to \mathcal{C}_N$ defined as follows.  Let $(x_0 < x_1)
\in N_1D_N^{\perp,\epsilon}$ with $x_0 = (a_0, G)$, $x_1 = (a_1, G)$.
Then let
\begin{align*}
  c(x_0 < x_1) = G|(a_0 - \epsilon, a_1 + \epsilon) \times \R^{N-1}.
\end{align*}
This defines a functor $D_N^{\perp,\epsilon} \to
\mathcal{C}_N^\epsilon$, and $c: D_N^\perp \to \mathcal{C}_N$ is
defined by taking the colimit.

The following lemma implies proposition~\ref{prop:3A1}.
\begin{lemma}\label{lemthm:equivalences}
  The induced maps
  \begin{align}
    \label{eq:8}
    Bi: BD^\perp_N & \to BD^\pitchfork_N\\
    \label{eq:9}
    Bu: BD^\pitchfork_N & \to D_N\\
    \label{eq:10}
    Bc: BD^\perp_N &\to B\mathcal{C}_N
  \end{align}
  are all weak equivalences.
\end{lemma}
\begin{proof}
  In fact~(\ref{eq:8}) and~(\ref{eq:10}) are both induced by
  degreewise weak homotopy equivalences on simplicial nerves.
  For~(\ref{eq:8}) this is obvious---straighten the morphisms near
  their ends.

  For~(\ref{eq:9}), notice that all maps
  \begin{align*}
    N_ku: N_k D_N^\pitchfork \to D_N
  \end{align*}
  are etale, and that for $G \in D_N$, the inverse image $u^{-1} (G)$
  is the set
  \begin{align*}
    \{a \in \R \mid G \pitchfork \{a\} \times \R^{N-1}\}
  \end{align*}
  which has the discrete topology, is non-empty, and totally ordered.
  It follows that $B(u^{-1}(G))$ is contractible, and the claim
  follows from lemma~\ref{lemma:Segal}.

  For~(\ref{eq:10}), suppose $P$ is a sphere and $f: P \to
  N_k\mathcal{C}_N$ a continuous map.  By compactness, $f$ maps into
  $N_k\mathcal{C}_N^\epsilon$ for some $\epsilon > 0$ so all graphs in
  the image of $f$ are elements of $\Phi((a_0 - \epsilon, a_k +
  \epsilon)\times \R^N)$.  Choose a diffeomorphism from $(a_0 -
  \epsilon, a_k + \epsilon)$ to $\R$ which is the identity on $(a_0 -
  \epsilon/2, a_k + \epsilon/2)$ and use that to lift $f$ to $P \to
  N_k D^\perp_N$.  We have constructed an inverse to $\pi_*(N_kc)$.
\end{proof}

We have proved the following result.
\begin{proposition}\label{prop:3A1}
  There is a weak equivalence
  \begin{align}
    B \mathcal{C}_N \simeq D_N.\tag*{\qed}
  \end{align}
\end{proposition}

A variation of the proof of proposition~\ref{prop:3A1} given above
will prove the following result.  The details are given below.

\begin{proposition}\label{prop:3A2}
  There is a weak equivalence
  \begin{align*}
    D_N \xrightarrow{\simeq} \Omega^{N-1}\Phi(\R^N).
  \end{align*}
\end{proposition}

The map $D_N \to \Omega^{N-1}\Phi(\R^N)$ is similar to the map
$\tau_N$ in definition~\ref{defn:alpha}.  First let $\R^{N-1} \times
D_N \to \Phi(\R^N)$ be given by the formula
\begin{align}\label{eq:11a}
  (v,G) \mapsto G - (0,v).
\end{align}
This extends uniquely to a continuous map $S^{N-1} \wedge D_N \to
\Phi(\R^N)$, and the adjoint of this map is the weak equivalence in
proposition~\ref{prop:3A2}.  This map is homotopic to a map
\begin{align}\label{eq:11b}
  D_N \to \Gamma((\R \times I^{N-1}, \R \times \partial
  I^{N-1}),\Phi^\mathrm{fib}(T\R^N)) \simeq \Omega^{N-1}\Phi(\R^N)
\end{align}
defined by ``scanning'', just like the map $\tau_N$ in
definition~\ref{defn:alpha} is homotopic to the map $\tilde \tau_N$ in
\eqref{eq:39}.

We give two proofs proposition~\ref{prop:3A2}.  The first is a direct
induction proof which is similar to the proofs of
propositions~\ref{prop:3A1} and \ref{prop:3B} above.  The second uses
Gromov's ``flexible sheaves'' \cite[section 2]{MR864505}.  While this
is somewhat heavy machinery, we believe it illuminates the relation
between scanning maps and Pontryagin-Thom collapse maps nicely.  For
the second proof, the crucial properties of $\Phi$ are the continuity
property expressed in proposition~\ref{prop:res-is-cont} and that
$\Phi$ is ``microflexible''.

\subsubsection{First proof}
\label{subsubsec:first-proof}

For $k = 0, 1, \dots, N$, let $D_{N,k} \subseteq \Phi(\R^N)$ be the
subspace
\begin{align*}
  D_{N,k} = \{ G \in \Phi(\R^N) \mid G \subseteq \R^k \times
  \Int(I^{N-k}) \},
\end{align*}
equipped with the subspace topology.  In particular $D_{N,0} = B_N$,
$D_{N,1} = D_N$ and $D_{N,N} = \Phi(\R^N)$.  The map $\R \times
D_{N,k-1} \to D_{N,k}$ given by
\begin{align*}
  (t,G) \mapsto G - (0,t,0)
\end{align*}
extends uniquely to a continuous map $S^1 \wedge D_{N,k-1} \to
D_{N,k}$ and we consider its adjoint
\begin{align}
  \label{eq:22}
  D_{N,k-1} \to \Omega D_{N,k}.
\end{align}
The composition of the maps~\eqref{eq:22} for $k=2,\dots, N$, is the
map $D_N \to \Omega^{N-1}\Phi(\R^N)$ of proposition~\ref{prop:3A2}.
\begin{proposition}\label{propthm:delooping-d-inductively}
  The map~\eqref{eq:22} is a weak equivalence for $k = 2, 3, \dots,
  N$.
\end{proposition}
Proposition~\ref{prop:3A2} then follows from
proposition~\ref{propthm:delooping-d-inductively} by induction.  The
proof of proposition~\ref{propthm:delooping-d-inductively} is given in
the lemmas~\ref{lemma:empty-intersections} and \ref{lem:Bc-Bi-Bu}
below.  The proofs of these lemmas are very similar to the proofs of
the propositions~\ref{prop:3B} and \ref{prop:3A1}, respectively.
\begin{definition}
  Let $k\geq 2$.
  \begin{enumerate}[(i)]
  \item Let $D^\pitchfork_{N,k}$ be the space of triples $(G,a,p) \in D_{N,k}
    \times \R^\delta \times (\R^\delta)^{k-1}$ satisfying
    \begin{align}\label{eq:34}
      \{p\} \times \{a\} \times \R^{N-k} \cap G = \emptyset.
    \end{align}
    Order $D^\pitchfork_{N,k}$ by declaring $(G,a,p) < (G',a',p')$ if
    and only if $G = G'$ and $a <a'$.
  \item Let $D^\perp_{N,k}\subseteq D^\pitchfork_{N,k}$ be the
    subspace, and sub poset, consisting of triples satisfying the
    further condition
    \begin{align}\label{eq:33}
      \R^{k-1} \times \{a\} \times \R^{N-k} \cap G = \emptyset.
    \end{align}
  \item Let $\mathcal{C}_{N,k}$ be the category whose space of objects
    is $\R^\delta \times (\R^\delta)^{k-1}$ and with morphism spaces
    given by
    \begin{multline*}
      \mathcal{C}_{N,k} ((a_0, p_0), (a_1, p_1)) = \\ \{ G \in \Phi(\R^N) \mid
      G \subseteq \R^{k-1} \times \Int([a_0,a_1] \times I^{N-k}) \},
    \end{multline*}
    when $a_0 \leq a_1$.  Composition is union of subsets.
  \item Let $i:D^\perp_{N,k} \to D^\pitchfork_{N,k}$ be the inclusion
    functor and let $u: D^\pitchfork_{N,k} \to D_{N,k}$ be the
    forgetful map.  Let $c: D^\perp_{N,k} \to \mathcal{C}_{N,k}$ be
    the functor given on morphisms by
    \begin{align*}
      c((G,a_0, p_0) \leq (G,a_1, p_1)) = G \cap \big(\R^{k-1} \times
      [a_0,a_1] \times \R^{N-k}\big).
    \end{align*}
  \end{enumerate}
\end{definition}

\begin{lemma}\label{lemma:empty-intersections}
  The maps
  \begin{align*}
    B\mathcal{C}_{N,k} \xleftarrow{Bc} BD^\perp_{N,k} \xrightarrow{ Bi }
    BD^\pitchfork_{N,k} \xrightarrow{ Bu } D_{N,k}
  \end{align*}
  are all weak equivalences.
\end{lemma}
\begin{proof}
  This is very similar to the proof of
  lemma~\ref{lemthm:equivalences}.  We first consider $Bu$.  For each
  $l$, $N_lu: N_lD^\pitchfork_{N,k} \to D_{N,k}$ is an etale map, and
  each fiber $u^{-1}(G)$ is a contractible poset (for each $a \in \R$,
  we can choose $p_a \in \R^{k-1}$ such that $(G,a,p_a) \in
  D^\pitchfork_{N,k}$.  The set of all $(G,a,p_a)$, $a \in \R$ forms a
  totally ordered cofinal subposet of $u^{-1}(G)$).  The result now
  follows from lemma~\ref{lemma:Segal}.

  For $Bi$ we claim that the inclusion $N_lD^\perp_{N,k} \to
  N_lD^\pitchfork_{N,k}$ is a deformation retract for each $l$.  A
  non-degenerate element $\chi \in N_l D^\pitchfork_{N,k}$ is given by
  an element $G \in D_{N,k}$, real numbers $a_0 < \dots < a_l$, and
  points $p_0, \dots, p_l \in \R^{k-1}$.  We will define a path from
  $\chi$ to a point in $N_lD^\perp_{N,k}$ depending continuously on
  $\chi$.  In essence, the construction is a parametrized version of
  the path constructed in lemma~\ref{lemma:Phi-is-connected}.

  For $r \in \R$, let $h_r: \R^{k-1} \to \R^{k-1}$ be the affine
  function given by
  \begin{align*}
    h_r(x) = x \prod_{i=0}^l(r-a_i)^2 + \sum_{i=0}^l p_i \prod_{j \neq
      i} \frac{r-a_j}{a_i - a_j}.
  \end{align*}
  Then $h_r$ is a diffeomorphism for $r \not\in \{a_0, \dots, a_l\}$
  and $h_{a_i}(x) = p_i$ for all $x$.  For $t \in [0,1]$, let $\phi_t:
  \R^{k-1} \times \R \times \R^{N-k} \to \R^{k-1} \times \R \times
  \R^{N-k}$ be the map given by
  \begin{align*}
    \phi_t(x,r,y) = ((1-t)x + t h_r(x),r,y)
  \end{align*}
  and let $G_t = (\phi_t)^*(G)$.  Then $G_0 = G$ and $G_1$ will
  satisfy~\eqref{eq:33} with respect to any $a = a_\nu$, $\nu \in
  \{0,1,\dots, l\}$.  This gives a continuous path in
  $N_lD^\pitchfork_{N,k}$,
  \begin{align*}
    t \mapsto \chi_t = (G_t, (a_0< \dots < a_l), (p_0, \dots, p_l)),
    \quad t \in [0,1],
  \end{align*}
  starting at $\chi_0 = \chi$ and ending at $\chi_1 \in N_l
  D^\perp_{N,k}$.

  Finally, $Bc$ is a homotopy equivalence because $N_lc$ is a homotopy
  equivalence for all $l$.  This is proved precisely as the analogous
  statement in lemma~\ref{lemthm:equivalences}.
\end{proof}

For the proof of lemma~\ref{lem:Bc-Bi-Bu}, recall that to a covariant
functor $F: C \to \mathrm{Spaces}$ there is an associated category $(C
\wr F)$.  When $F$ is contravariant we let $(F \wr C) = (C^\mathrm{op}
\wr F)^\mathrm{op}$.

\begin{lemma}\label{lem:Bc-Bi-Bu}
  Let $k\geq 2$.
  \begin{enumerate}[(i)]
  \item\label{item:1} $\mathcal{C}_{N,k}((a_0,p_0), (a_1,p_1)) \cong D_{N,k-1}$
    whenever $a_0 < a_1$.
  \item\label{item:2} Composition with any morphism $G:
    (a_1,p_1) \to (a_2, p_2)$ induces a homotopy equivalence
    \begin{align*}
      \mathcal{C}_{N,k}((a_0,p_0), (a_1,p_1)) \xrightarrow{G\circ}
      \mathcal{C}_{N,k}((a_0,p_0), (a_2,p_2)).
    \end{align*}
    Similarly for composing from the right.
  \item\label{item:3} There is a weak equivalence
    \begin{align*}
      D_{N,k} \simeq \Omega B \mathcal{C}_{N,k-1}.
    \end{align*}
  \end{enumerate}
\end{lemma}
\begin{proof}
  \emph{(\ref{item:1})} There is an obvious homeomorphism that
  stretches the interval $[a_0,a_1]$ to $[0,1]$.

  \emph{(\ref{item:2})} is clear in the case $G = \emptyset$.  The
  space $D_N = D_{N,1} \simeq B\mathcal{C}_N$ is connected (given any
  two objects, there is a morphism between them).  Similarly $D_{N,k}
  \simeq B \mathcal{C}_{N,k}$ is connected for $k \geq 2$.  Therefore
  composition with any $G \in \mathcal{C}_{N,k}((a_1,p_1),
  (a_2,p_2))\simeq D_{N,k-1}$ is homotopic to composition with $G =
  \emptyset$.

  \emph{(\ref{item:3})}.  Consider for each $n \in \Z$ the object
  $(n,0)$ in $\mathcal{C}_{N,k}$, and let $F_n: \mathcal{C}_{N,k} \to
  \mathrm{Spaces}$ be the functor
  \begin{align*}
    F_n = \mathcal{C}_{N,k}( - , (n,0)).
  \end{align*}
  The morphism $\emptyset: (n,0) \to (n+1,0)$ induces a natural
  transformation $F_n \to F_{n+1}$, and we let
  \begin{align*}
    F_\infty ((a,p)) = \hocolim_{n \to \infty} F_n((a,p)).
  \end{align*}
  Then $\mathrm{id}_{(n,0)}$ is a final object of the category $(F_n
  \wr \mathcal{C}_{N,k})$, so $B(F_n \wr \mathcal{C}_{N,k})$ and
  $B(F_\infty \wr \mathcal{C}_{N,k})$ are contractible.

  From \emph{(\ref{item:1})} and \emph{(\ref{item:2})} we get that
  $F_\infty((a,p)) \simeq D_{N,k-1}$, and that any morphism in
  $\mathcal{C}_{N,k}$ induces a homotopy equivalence $F_\infty((a_1,
  p)) \to F_\infty((a_0, p))$.  Then the simplicial map
  $N_\bullet(F_\infty \wr \mathcal{C}_{N,k}) \to
  N_\bullet(\mathcal{C}_{N,k})$ satisfies the hypothesis of
  \cite[1.6]{MR0353298}, so the geometric realization
  \begin{align*}
    B(F_\infty \wr \mathcal{C}_{N,k}) \to B\mathcal{C}_{N,k}
  \end{align*}
  is a quasifibration.  Thus for any $(a,p) \in N_0\mathcal{C}_{N,k}$,
  the inclusion of an actual fiber over $(a,p)$ into the homotopy
  fiber is a weak equivalence.  The actual fiber over $(0,0)$ is
  $F_\infty((0,0)) \simeq D_{N,k-1}$ and the homotopy fiber is
  equivalent to $\Omega_{(0,0)} B\mathcal{C}_{N,k}$.
\end{proof}

\subsubsection{Second proof}\label{subsubsec:second-proof}
The sheaf $\Phi$ is an example of an \emph{equivariant, continuous
  sheaf} in the terminology of \cite{MR864505}.  This means that
$\Phi$ is continuously functorial with respect to embeddings (not just
inclusions) of open subsets of $\R^N$, cf.\
proposition~\ref{prop:res-is-cont}.  In particular, $\Diff(U)$ acts
continuously on $\Phi(U)$.  To such a sheaf on a manifold $V$ there is
an associated sheaf $\Phi^*$ and a map of sheaves $\Phi \to \Phi^*$.
Up to homotopy, $\Phi^*(V)$ is the space of global sections of the
fiber bundle $\Phi^\mathrm{fib}(TV)$ defined in
section~\ref{subsec:param-pontrj-thom}, and the inclusion
\begin{align}\label{eq:29-Gromov}
  \Phi(V) \to \Phi^*(V) \simeq \Gamma(V,\Phi^\mathrm{fib}(TV))
\end{align}
is a scanning map induced by an ``exponential'' map on $V$, similar to
the map~\eqref{eq:39}.  Gromov, in \cite[section 2.2.2]{MR864505},
proves that \eqref{eq:29-Gromov} is a weak homotopy equivalence when
$V$ is \emph{open}, i.e.\ all connected components are non-compact,
and $\Phi$ is \emph{microflexible} (we recall the definition below).
This also holds in a relative setting $(V,\partial V)$.  In particular
we can use $(V, \partial V) = (\R \times I^{N-1}, \R \times \partial
I^{N-1})$, in which case~\eqref{eq:29-Gromov} specializes
to~\eqref{eq:11b}.

That the sheaf $\Phi$ is microflexible means that for each inclusion
of compact subsets $K' \subseteq K \subseteq \R^N$, each open $U,
U'$ with $K' \subseteq U' \subseteq U \supseteq K$, and each diagram
\begin{align}\label{eq:30} 
  \xymatrix{
    {P \times\{0\}} \ar[r]^-{h}\ar[d] & {\Phi(U')}\ar[d] \\
    {P \times [0,1]} \ar[r]^-{f} & {\Phi(U)}
  }
\end{align}
with $P$ a compact polyhedron, there exists an $\epsilon > 0$ and an
initial lift $P \times [0,\epsilon] \to \Phi(U')$ of $f$ extending
$h$, after possibly shrinking $U\supseteq K$ and $U'\supseteq K'$.

In this subsection we prove that the sheaf of graphs is microflexible.
Then Gromov's $h$-principle implies that the map~(\ref{eq:11b}) above
is an equivalence for all $N$.

\begin{proposition}\label{prop:microflex}
  Let $K \subseteq U$ be compact and $P$ a polyhedron.  Let $f: P
  \times [0,1] \to \Phi(U)$ be continuous.  Then there exists an
  $\epsilon > 0$ and a continuous map $g: P \times [0,\epsilon] \to
  \Phi(U)$ with the following properties.
  \begin{enumerate}[(i)]
  \item\label{item:16} The map $f|P \times [0,\epsilon]$ agrees with
    $g$ near $K$,
  \item\label{item:17} the map $g|P \times \{0\}$ agrees with $f|P
    \times \{0\}$,
  \item\label{item:18} there exists a compact subset $C\subseteq U$
    such that the map
    \begin{align}\label{eq:31}
      P \times[0,\epsilon] \stackrel{g}{\longrightarrow} \Phi(U)
      \stackrel{\mathrm{res}}{\longrightarrow} \Phi(U - C)
    \end{align}
    factors through the projection $\mathrm{pr}: P\times [0,\epsilon]
    \to P$.
  \end{enumerate}
\end{proposition}
Proposition~\ref{prop:microflex} immediately implies microflexibility.
Indeed, given maps as in diagram~\eqref{eq:30}, the composition
\begin{align*}
  P \times [0,\epsilon] \xrightarrow{\mathrm{pr}} P \times \{0\}
  \xrightarrow{ h} \Phi(U') \to \Phi(U' - C)
\end{align*}
will agree with $g: P\times [0,\epsilon] \to \Phi(U)$ on the overlap
$U \cap (U' - C) = U-C$, so they can be glued together to a map $P
\times [0,\epsilon] \to \Phi(U')$.  The glued map is the initial lift
in diagram~\eqref{eq:30}.

\begin{proof}[Proof for $P$ a point] We are given a continuous path
  $f: [0,1] \to \Phi(U)$.  Let $C\subseteq U$ be compact with $K
  \subseteq \Int(C)$ and choose $\tilde\tau: U \to [0,1]$ with
  $\tilde\tau = 1$ near $K$ and with $\supp(\tilde\tau) \subseteq C$
  compact.  For each of the finitely many vertices $q \in \VV(f(0))
  \cap (\supp(\tilde\tau) - K)$, choose a function $\rho_q: U \to
  [0,1]$ which is 1 near $q$, such that the sets $\supp(\rho_q)$ have
  compact support in $U - K$ and are mutually disjoint.  Let $\tau: U
  \to [0,1]$ be the function
  \begin{align*}
    \tau(v) = \tilde\tau (v) + \sum_q \rho_q(v)\big(\tilde\tau(v) -
    \tilde\tau(q)\big).
  \end{align*}
  Then $\tau: U \to [0,1]$ is locally constant outside a compact
  subset of $U - (K \cup \VV(f(0)))$.

  Continuity of $f$ gives a graph epimorphism $\phi_t: f(t)
  \dashrightarrow f(0)$ for $t$ sufficiently close to $0$, defined and
  canonical near $C$.  Let $g(t)$ be the image of the map
  \begin{align*}
    f(t) & \to U\\
    x & \mapsto \tau(x)x + (1-\tau(x))\phi(x).
  \end{align*}
  For $t\in[0,1]$ sufficiently close to 0, this defines an element
  $g(t) \in \Phi(U)$ satisfying (\ref{item:16}), (\ref{item:17}),
  (\ref{item:18}).
\end{proof}
\begin{proof}[General case] To make the above argument work in the
  general case (parametrized by a compact polyhedron $P$), we need
  only explain how to choose the function $\tau: P \times U \to
  [0,1]$.  For each $p\in P$, the above construction provides a
  $\tau_p: U \to [0,1]$ that works for $f|\{p\} \times [0,1]$ (i.e.\
  $\tau_p(x,u)$ is independent of $u$ near vertices of $f(p,0)$).  The
  same $\tau_p$ will work for $f|\{q\} \times [0,1]$ for all $q$ in a
  neighborhood $W_p \subseteq P$ of $p$.  Choose a partition of unity
  $\lambda_p: P \to [0,1]$ subordinate to the open covering by the
  $W_p$.  Then let
  \begin{align}\tag*{\qedhere}
    \tau(q,v) = \sum_p \lambda_p(q) \tau_p(v).
  \end{align}
\end{proof}

\subsection{The positive boundary subcategory}
\label{subsec:posit-bound-subc}

The condition on morphisms in the positive boundary subcategory
$\mathcal{C}_N^\partial \subseteq \mathcal{C}_N$
(definition~\ref{defn:pos-bd}) ensures that any morphism $G: (a_0,
A_0, \lambda_0) \to (a_1, A_1, \lambda_1)$ is connected when
$|A_1|=1$.  This will allow us to use homological stability to prove
the ``group completion'' result in proposition~\ref{prop:3B} using
\cite{MR0402733}, much as it was done in the parallel case of
two-dimensional manifolds, \cite{MR1474157}.
\begin{lemma}\label{lemthm:pos-bd-subcat}
  Let $c_0 = (a_0,A_0, \lambda_0)$ and $c_1 = (a_1, A_1, \lambda_1)$
  be two objects of $\mathcal{C}_\infty^\partial$, with $a_0 < a_1$ and
  $|A_1| = 1$.  Then
  \begin{align*}
    \mathcal{C}_\infty^\partial(c_0,c_1) \simeq \coprod_n
    BA_n^{1+|A_0|}.
  \end{align*}
\end{lemma}
\begin{proof}
  The surjectivity of $A_1 \to \pi_0(G)$ implies that $G$ is
  connected.  Then the lemma follows from
  theorem~\ref{thmcor:main-decomposition}.
\end{proof}

\begin{proposition}
  \label{prop:3B}
  There is a homology equivalence
  \begin{align*}
    \Z \times B \Aut_\infty \to \Omega B \mathcal{C}_\infty^\partial.
  \end{align*}
\end{proposition}
\begin{proof}
  This is very similar to lemma~\ref{lem:Bc-Bi-Bu}(\ref{item:3}), but
  with homology equivalences instead of weak homotopy equivalences.
  We sketch the proof.  See \cite[chapter 7]{math.AT/0605249} for more
  details.

  Let $A \subseteq \Int(I^N)$ be a one-point set.  For $n \in
  \N\subseteq \R$, let
  \begin{align*}
    F_n: (\mathcal{C}_\infty^\partial)^\mathrm{op} \to \mathrm{Spaces}
  \end{align*}
  be the functor $F_n = \mathcal{C}_\infty^\partial( - ,(n,A,0))$.
  $(F_n \wr \mathcal{C}_\infty^\partial)$ has $\mathrm{id}_{(n,A,0)}$
  as final object, so $B(F_n \wr \mathcal{C}_\infty^\partial)$ is
  contractible.  Choose morphisms $G_n: (n,A,0) \to (n+1,A,0)$ in
  $\mathcal{C}_\infty^\partial$ with first Betti number $b_1(G_n) =
  1$.  This defines a direct system
  \begin{align*}
    F_1 \xrightarrow{G_1} F_2 \xrightarrow{G_2} \dots
    \xrightarrow{G_n} F_n \xrightarrow{G_{n+1}} \dots
  \end{align*}
  and we let $F_\infty(x) = \hocolim_n F_n(x)$.  Then $B(F_\infty \wr
  \mathcal{C}_\infty^\partial) = \hocolim_n B(F_n \wr
  \mathcal{C}_\infty^\partial)$ is still contractible.

  Lemma~\ref{lemthm:pos-bd-subcat} gives a homotopy equivalence for
  each object $c = (a,A_0, \lambda)$
  \begin{align*}
    F_\infty(c) \simeq \Z \times B A_\infty^{1+|A_0|}
  \end{align*}
  and by theorem~\ref{thm:HVW}, the functor $F_\infty:
  (\mathcal{C}_\infty^\partial)^\mathrm{op} \to \mathrm{Spaces}$ maps
  every morphism to a homology equivalence.  This implies that the
  simplicial map $N_\bullet(F_\infty\wr \mathcal{C}_\infty^\partial)
  \to N_\bullet \mathcal{C}_\infty^\partial$ satisfies the assumption
  of \cite[proposition 4]{MR0402733} and therefore that the geometric
  realization
  \begin{align}\label{eq:21}
    B(F_\infty \wr \mathcal{C}_\infty^\partial) \to
    B\mathcal{C}_\infty^\partial
  \end{align}
  is a homology fibration in the sense of \cite{MR0402733}.  Thus for
  any $(a,A,\lambda) \in N_0\mathcal{C}_\infty^\partial$, the
  inclusion of an actual fiber over $(a,A)$ of~\eqref{eq:21} into the
  homotopy fiber is a homology equivalence.  The actual
  fiber over $(0,\emptyset,0)$ is $F_\infty((0,\emptyset,0)) \simeq \Z
  \times B \Aut_\infty$.  Since $B(F_\infty \wr
  \mathcal{C}_\infty^\partial)$ is contractible, the homotopy fiber is
  equivalent to $\Omega_{(0,\emptyset,0)}B
  \mathcal{C}_\infty^\partial$.
\end{proof}

The following proposition is proved in several steps.  The proof
occupies the rest of this section, and is similar to \cite[chapter
6]{math.AT/0605249}.

\begin{proposition}\label{prop:3C}
  The inclusion induces a weak equivalence
  \begin{align*}
    B\mathcal{C}_\infty^\partial \xrightarrow{\simeq} B \mathcal{C}_\infty.
  \end{align*}
\end{proposition}

For $G \in D_N$, we shall write $f_G: G\to \R$, or just $f$, for the
restriction to $G$ of the projection $\R \times \Int(I^{N-1}) \to \R$.
\begin{definition}
  Let $G \in D_N$ and $p \in f^{-1}((-\infty,0])$.  Define $f^-(p)\in
  [-\infty,f(p)]$ as
  \begin{align*}
    f^-(p) = \max_\gamma \min_{t\in[0,1]} f\gamma(t)
  \end{align*}
  where the maximum is taken over paths $\gamma: [0,1] \to G$
  satisfying $\gamma(0) = p$ and $f\gamma(1)>0$.  We let $f^-(p) = -
  \infty$ if no such $\gamma$ exists.  Let $A_G\subseteq
  f^{-1}((-\infty,0])$ be the closure of the set of points for which
  $f^-(p) < f(p)$.  Let $B_G\subseteq G$ be the union of $A_G$, the
  set of vertices, and the set of edgewise critical points of $f$.
  Let $R_G = (-\infty,0] - f(B_G)$ and
  \begin{align*}
    D_N^\partial = \{G \in D'_N \mid R_G \neq \emptyset\}.
  \end{align*}
\end{definition}
For fixed $r \in (-\infty,0]$, the set $\{G \in D_N\mid r \in R_G\}$
is an open subset of $D_N$.

\begin{lemma}\label{lemprop:pos-bdy}
  There is a weak equivalence $D_N^\partial \simeq
  B\mathcal{C}_N^\partial$.
\end{lemma}

\begin{proof}
  This is completely analogous to the proof of
  theorem~\ref{lemthm:equivalences} in
  section~\ref{subsec:poset-model-graph}.  It uses the subposet
  $D_N^{\partial,\pitchfork}$ of $D_N^\pitchfork$ consisting of
  $(a,G)$ with $G \in D_N^\partial$ and $a \in R_G$ and the poset
  $D_N^{\partial, \perp} = D_N^\perp \cap D_N^{\partial,\pitchfork}$.
  As in the proof of theorem~\ref{lemthm:equivalences} we have levelwise
  equivalences
  \begin{align*}
    N_\bullet D_N^{\partial, \perp} \xrightarrow{\simeq} N_\bullet
    D_N^{\partial, \pitchfork}, \quad N_\bullet D_N^{\partial, \perp}
    \xrightarrow{\simeq} N_\bullet \mathcal{C}_N^\partial,
  \end{align*}
  and the equivalence $BD_N^{\partial, \pitchfork} \to D_N^\partial$
  uses lemma~\ref{lemma:Segal}.
\end{proof}
Proving proposition~\ref{prop:3C} now amounts to the inclusion
$D_N^\partial \subset D_N$ being a weak equivalence.  This is done in
the lemmas~\ref{lemproposition:no-compact-components} and
\ref{lemprop:get-rid-of-local-max} below. 
\begin{definition}
  Let $D_N'\subseteq D_N$ be the subset consisting of graphs $G$ for
  which no path component of $G$ is compact.
\end{definition}

\begin{lemma}\label{lemproposition:no-compact-components}
  The inclusion $D'_N \to D_N$ is a weak equivalence.
\end{lemma}
\begin{proof}
  For a given $G \in D'_N$, we can assume, after possibly perturbing
  the function $f$ a little, that no connected component of $f$ is
  contained in $f^{-1}(0)$.  Then we can choose an $\epsilon > 0$
  small enough that no connected component of $f^{-1} ((-\epsilon,
  \epsilon)) \subseteq G$ is compact.  For $t \in [0,1]$ let $h_t: \R
  \to \R$ be an isotopy of embeddings with $h_0 = \id$ and $h_1(\R) =
  (-\epsilon,\epsilon)$.  Let $H_t = h_t \times \id: \R \times
  \R^{N-1} \to \R \times \R^{N-1}$.  Then
  \begin{align*}
    t \mapsto G_t = H_t^*(G)
  \end{align*}
  defines a continuous path $[0,1] \to D_N$, starting at $G_0 = G$ and
  ending in $G_1 \in D_N'$.

  This proves that the relative homotopy group $\pi_k(D_N, D_N')$ is
  trivial for $k = 0$.  The case $k >0$ is similar: Given a continuous
  map of pairs
  \begin{align*}
    q: (\Delta^k, \partial \Delta^k) \to (D_N, D_N')
  \end{align*}
  we can first perturb $q$ a little, such that for all $x \in
  \Delta^k$, no connected component of $q(x)$ is contained in
  $f^{-1}(0)$, and then stretch a small interval $(-\epsilon,
  \epsilon)$.
\end{proof}

The relevance of the condition that no connected component of $G \in
D_N'$ be compact lies in the following definition.
\begin{definition}
  For $G \in D_N$ let $\hat G = G \amalg \{+\infty, -\infty\}$.  Then
  $f$ extends to $f: \hat G \to [-\infty,\infty]$, and we equip $\hat
  G$ with the coarsest topology in which $G\subseteq \hat G$ has the
  subspace topology and $f: \hat G \to [-\infty, \infty]$ is
  continuous.  (In other words, a sequence of points $x_n \in G$,
  $n\in\N$, converges to $\pm \infty \in \hat G$ if and only if
  $f(x_n) \to \pm \infty$.)  An \emph{escape to $+\infty$} is a path
  $\gamma: [0,1] \to \hat G$ such that $\gamma(0) = p$ and $\gamma(1)
  = +\infty$.  An escape to $-\infty$ is defined similarly.
\end{definition}
Given $G$ and $p$, an escape to either $+\infty$ or $-\infty$ exists
if and only if the path component of $G$ containing $p$ is
non-compact.  Let us also point out that a path $\gamma: [0,1] \to
\hat G$ is uniquely given by its restriction $[0,1] \dashrightarrow
G$, defined on $\gamma^{-1}(G)\subseteq [0,1]$.

\begin{remark}
  \label{remark:local-paths}
  The statement of lemma~\ref{lemproposition:no-compact-components} is
  that any map of pairs $q:(\Delta^k, \partial \Delta^k) \to (D_N,
  D_N')$ is homotopic to a map $q'$ such that for any $x \in \Delta^k$
  there exists an escape to $\pm \infty$ from $p \in q'(x)$.  In fact,
  essentially the same proof gives a slightly stronger statement,
  namely that such escapes exist locally in $\Delta^k$ (not just
  pointwise).

  Indeed, if $p \in f^{-1}((-\epsilon,\epsilon))$ and $\gamma: [0,1]
  \to G$ is a path with $\gamma(0) = p$ and $|f\gamma(1)| >\epsilon$
  and $H_t$ is the isotopy from the proof of
  lemma~\ref{lemproposition:no-compact-components}, then $H_1^{-1}
  \circ \gamma: [0,1] \dashrightarrow G_1 = H_1^*(G)$ is an escape
  from $H_1^{-1}(p)$ to either $+\infty$ or to $-\infty$.  If $G =
  q(x_0)$ for some $x_0 \in \Delta^k$, then the path $\gamma$ can be
  extended locally to $\Gamma: U_x \times [0,1] \to \R^N$ for a
  neighborhood $U_x\subseteq \Delta^k$ of $x$, such that $\Gamma(x,t)
  \in q(x)$ and $\Gamma(x_0, -) = \gamma$.  Then $H_1^{-1} \circ
  \Gamma$ is a family of escapes to $+\infty$ or $- \infty$, defined
  locally near $x_0$.
\end{remark}

\begin{lemma}\label{lemprop:get-rid-of-local-max}
  The inclusion $D^\partial_\infty \to D_\infty'$ is a weak homotopy
  equivalence.
\end{lemma}
\begin{proof}
  We prove that for $k\geq 0$, any map of pairs
  \begin{align}
    \label{eq:23}
    q: (\Delta^k, \partial \Delta^k) \to (D_\infty', D_\infty^\partial)
  \end{align}
  is homotopic to a map into $D_\infty^\partial$.

  Consider first the case $k=0$.  Let $G = q(1)$.  Choose $a,b \in \R$
  with $a<0<b$ and $G \pitchfork \{a,b\} \times \R^\infty$.  If $G$
  satisfies the condition that
  \begin{align}\label{eq:24}
    \text{$\pi_0(f^{-1}(b)) \to \pi_0(f^{-1}([a,b]))$ is surjective}
  \end{align}
  then $[a, a + \epsilon] \subseteq R_G$ for some $\epsilon > 0$, and
  hence $G\in D_\infty^\partial$.  For general $G \in D_\infty'$ we
  will construct a path $h: [0,1] \to D_\infty'$ with $h(0) = G$ and
  such that $h(1)$ satisfies~\eqref{eq:24}.

  \begin{figure}[htbp]
    \centering
    \includegraphics{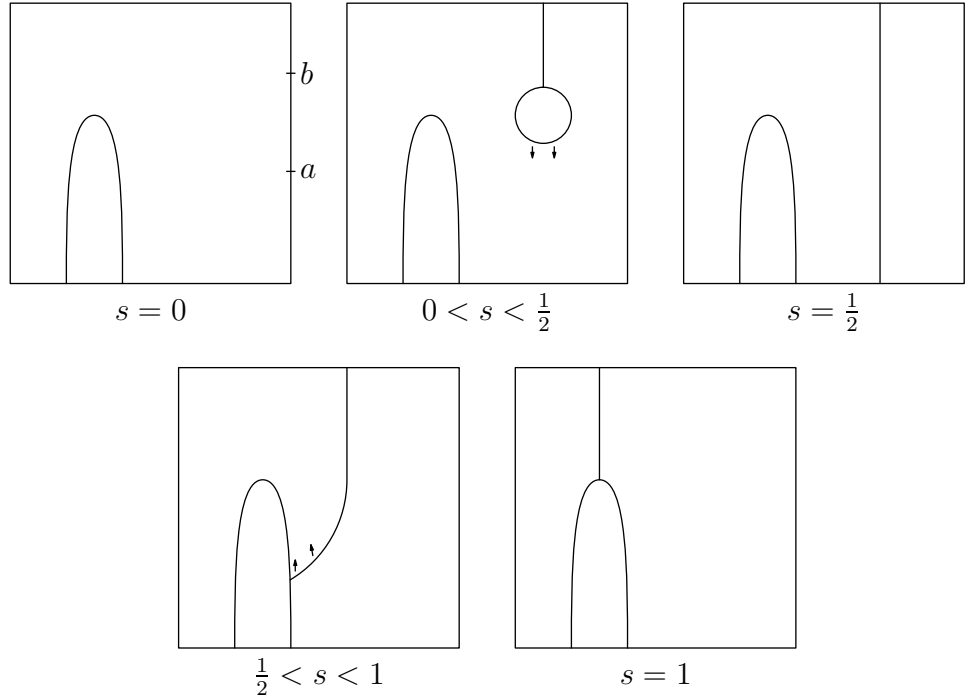}
    \caption{$h(s)$ for various $s \in [0,1]$}
    \label{fig:1}
  \end{figure}
  Let $p \in f^{-1}([a, b])$, and let $\gamma: [0,1] \dashrightarrow
  G$ be an escape from $p$ to $-\infty$.  A typical such $G$ is
  depicted in the first picture in the cartoon in figure~\ref{fig:1}
  which also depicts a path $h = h_\gamma: [0,1] \to D_\infty'$.  The
  pictures show part of the graph $h(s) \in D'_\infty$ for various $s
  \in [0,1]$.  $f_{h(s)}: h(s) \to \R$ is the height function
  (projection onto vertical axis) in the pictures.  At time $s = 0$
  the graph $G = h(0)$ has a local maximum in $f^{-1}([a,b])$.  At
  times $s \in (0,\frac12)$, the graph $h(s)$ is the disjoint union of
  $G$ and a ``gallow''.  At times $s \in [\frac12, 1]$, $h(s)$ is
  obtained from $G$ by attaching an extra edge at the point
  $\gamma(2-2s) \in \hat G$.  The path $h$ depends on two choices.
  Most importantly, it depends on the escape $\gamma$ from $p$ to
  $-\infty$ along which to ``slide'' the attached extra edge.
  Secondly, we have only defined the graph $h(s)$ abstractly; to get
  an element of $D_\infty'$, we choose an embedding $h(s) \subset \R
  \times \Int(I^{N-1})$ extending the inclusion of $G\subset h(s)$.
  Such an embedding always exists when $N = \infty$, so we suppress it
  from the notation.

  The path $h_\gamma$ has two convenient properties.  Firstly we have
  the inclusion $G \subseteq h(s)$ for all $s \in [0,1]$.  As
  constructed in the cartoon in figure~\ref{fig:1}, all local maxima
  of $f_{h(s)}: h(s) \to \R$ are in $G\subseteq h(s)$, so the subset
  \begin{align*}
    R_G \cap R_{h(s)} \subseteq R_G
  \end{align*}
  is open and dense.  In particular $R_{h(s)}$ is non-empty if $R_G$
  is non-empty.  So the path $h: [0,1] \to D_\infty'$ runs entirely in
  $D_\infty^\partial$, provided $h(0) \in D_\infty^\partial$.
  Secondly, at time $s = 1$, the graph $G_1 = h(1)$ is obtained from
  $G$ by attaching an extra edge extending from $p \in G$ to
  $+\infty$.  This assures that if $x \in f^{-1}([a, b])$ is in
  the same path component as $p$, then there is an escape from $x$ to
  $+\infty$ which stays inside $f_{G_1}^{-1}([a, \infty))$.  In
  particular $f_{G_1}^-(x) \geq a$.  If $x$ is sufficiently close to
  $f^{-1}(a)$, then by transversality we will even have
  $f_{G_1}^-(x) = x$.

  There is a similar construction if $\gamma: [0,1] \to\hat G$ is an
  escape from $p$ to $+\infty$, only easier: Let $h_\gamma(s)$ be the
  graph obtained from $G$ by attaching an extra edge extending to
  $+\infty$ at the point $\gamma(1-s) \in \hat G$.

  The same construction can be applied to attach several extra edges
  at the same time.  Let $X \subseteq f_G^{-1}((a,b))$ be a finite
  subset and let $\Gamma: X \times [0,1] \to \hat G$ be such that
  $\Gamma(p,-)$ is an escape from $p$ to $\pm \infty$.  Then the above
  construction gives a path $h = h_\Gamma: [0,1] \to D_\infty'$ such
  that $h_\Gamma(1)$ is obtained from $G = h_\Gamma(0)$ by attaching
  an extra edge extending to $+\infty$ at all the points $p \in X$.
  If $X\subseteq f^{-1}((a,b))$ is chosen such that the inclusion
  \begin{align*}
    X \amalg f^{-1}(b) \to f^{-1}([a,b])
  \end{align*}
  induces a surjection in $\pi_0$, then the resulting graph $G_1 =
  h(1)$ will satisfy $f_{G_1}^-(x) = x$ for all $x$ in a sufficiently
  small neighborhood of $f^{-1}(a)$, and hence $R_{G_1} \neq
  \emptyset$, so $h(1) \in D_\infty^\partial$.  This finishes the
  proof for $k=0$.

  For $k>0$ we will give a parametrized version of the above argument.
  If $\Gamma: X \times [0,1] \to \hat G$ and $\Gamma': X' \times [0,1]
  \to \hat G$ are two sets of escapes to $\pm \infty$, then we get two
  paths $h = h_\Gamma: [0,1] \to D_\infty'$ and $h' = h_{\Gamma'}:
  [0,1] \to D_\infty'$.  We will say that $h$ and $h'$ are
  \emph{compatible} if the inclusion
  \begin{align*}
    G \subseteq h(s) \cap h'(s')
  \end{align*}
  is an equality for all $s,s' \in [0,1]$.  In this case we will have
  that $h(s) \cup h'(s')$ is the pushout of the diagram $h'(s')
  \leftarrow G \rightarrow h(s)$ and that $h(s)\cup h'(s') \in
  D_\infty'$.  This produces a continuous map
  \begin{align*}
    [0,1]^2 & \to D_\infty'\\
    (s,s') & \mapsto h(s) \cup h'(s').
  \end{align*}
  More generally a finite set of escapes $\Gamma_j: X_j \times [0,1]
  \to \hat G$, $j \in J$ produces a finite set of paths $h_j: [0,1]
  \to D_\infty'$, and their union gives a map
  \begin{align*}
    h_J: [0,1]^J &\to D_\infty'\\
    (s_j)_{j \in J} & \mapsto \cup_j h(s_j),
  \end{align*}
  provided the $h_j$ are \emph{pairwise compatible}, i.e.\ that $h_i$
  and $h_j$ are compatible for all $i,j \in J$ with $i \neq j$.

  We will have $h_J(t) \in D_\infty^\partial$ as long as \emph{at
    least} one of the coordinates of $t \in [0,1]^J$ is 1.  This is an
  important property of the homotopies constructed from the cartoon:
  If we have already attached enough extra edges to ensure $G \in
  D_\infty^\partial$, then attaching even more edges will not destroy
  the property of being in $D_\infty^\partial$ (even if we stop in the
  middle of the attaching process).

  Now let $k\geq 1$, and let $q$ be as in~\eqref{eq:23}.  We will use
  the above construction to prove that $q$ is homotopic to a map into
  $D_\infty^\partial$, and hence that
  $\pi_k(D'_\infty,D_\infty^\partial)$ vanishes.  If $x \in \Delta^k$
  has $q(x) \not\in D_\infty^\partial$, then the proof in the case
  $k=0$ gives a path $h = h_\Gamma$ from $q(x)$ to a point in
  $D_\infty^\partial$, depending on a family $\Gamma$ of escapes to
  $\pm \infty$.  Extending $\Gamma$ to a continuous family $\Gamma_y :
  X \times [0,1] \to \widehat{q(y)}$, $y \in U_x$, parametrized by a
  neighborhood $U_x$ of $x$, we get a homotopy
  \begin{align*}
    h_x: U_x \times [0,1] \to D_\infty'
  \end{align*}
  starting at $q|U_x$ and ending in a map $U_x \to D_\infty^\partial$.
  The extension of $\Gamma$ to a continuous family $\Gamma_y$, $y \in
  U_x$, can be assumed to exist by remark~\ref{remark:local-paths}.
  Thus we get an open covering of $\Delta^k$ by the sets $U_x$, $x \in
  \Delta^k$, and corresponding homotopies $h_x$.  We now explain how
  to glue all these together.

  Choose a triangulation $K$ of $\Delta^k$ so fine that for all $v \in
  \Ver(K)$ we have $\st(v) \subseteq U_x$ for some $x \in \Delta^k$.
  Let $\sd K$ denote the barycentric subdivision of $K$.  For $\sigma
  \in \Ver(\sd K) = \Simp(K)$ we write $\st(\sigma)$ for the open star
  of $\sigma$ as a 0-simplex of $\sd(K)$.  We write $\dim(\sigma)$ for
  the dimension as a simplex of $K$, and $\tau < \sigma$ if $\tau$ is
  a proper face (in $K$) of $\sigma$.  Choose bump functions
  $\lambda_\sigma: |K| \to [0,1]$, $\sigma \in \Ver(\sd K)$ with the
  properties
  \begin{enumerate}[(i)]
  \item $\supp(\lambda_\sigma) \subseteq \st(\sigma)$ for all
    $\sigma$,
  \item $|K| = \cup_{\sigma} \Int( \lambda_\sigma^{-1}(1))$.
  \end{enumerate}
  Then we have $\supp(\lambda_\sigma) \cap \supp(\lambda_\tau) =
  \emptyset$ unless $\tau < \sigma$ (or $\sigma = \tau$ or $\sigma <
  \tau$).  For a simplex $\chi = (\sigma_0 < \sigma_1 < \dots <
  \sigma_l) \in \Simp(\sd K)$, we have a corresponding geometric
  simplex
  \begin{align*}
    |\chi| \subseteq |\sd K| = |K|
  \end{align*}
  and $\lambda_\tau$ vanishes on this subspace unless $\tau$ is a
  vertex of $\chi$.

  For each $\sigma \in \Ver(\sd K)$, choose an $x \in \Delta^k$ with
  $\st(\sigma) \subseteq U_x$.  Let $h_\sigma = h_x | (\st(\sigma)
  \times [0,1])$.  This gives a homotopy
  \begin{align*}
    h_\sigma: \st(\sigma) \times [0,1] \to D_\infty',
  \end{align*}
  starting at $g|\st(\sigma)$ and ending in a map $\st(\sigma) \to
  D_\infty^\partial$.  Proceeding by induction on $\dim(\sigma)$, we
  can assume that $h_\sigma$ is compatible with $h_\tau$ over
  $\supp(\lambda_\sigma) \cap \supp(\lambda_\tau)$ for all faces $\tau
  < \sigma$ (we can assume this because $N = \infty$).  Then define
  $H_\chi$ as the composition
  \begin{align*}
    |\chi| \times [0,1] \to |\chi| \times [0,1]^{l+1} \to D_\infty',
  \end{align*}
  where the first map is given by
  \begin{align*}
    (x,s) \mapsto (x,(\lambda_{\sigma_0}(x), \dots, \lambda_{\sigma_l}(x))s)
  \end{align*}
  and the second is given as $\cup_{i=0}^l h_{\sigma_i}$.  This is
  well-defined because the $\sigma_i$'s are proper faces of each
  other, so the homotopies $h_{\sigma_i}$ are compatible.  The
  homotopies $H_\chi$ glue together to a homotopy
  \begin{align*}
    H: |K| \times [0,1] \to D_\infty',
  \end{align*}
  starting at $q$ and ending in a map $H(-,1): \Delta^k \to
  D_\infty^\partial$.
\end{proof}

\section{Homotopy type of the graph spectrum}
\label{sec:homotopy-type-Phi}

The main result in this section is the following, which will finish
the proof of theorem~\ref{thm:MW}.
\begin{theorem}\label{thm:homotopy-type-Phi}
  We have an equivalence of spectra $\mathbf{\Phi} \simeq S^0$ and
  hence a weak equivalence
  \begin{align*}
    \Omega^\infty\mathbf{\Phi} \simeq QS^0.
  \end{align*}
\end{theorem}

Let us first give an informal version of the proof.  Since any
$\epsilon$-neighborhood of $0 \in \R^N$ can be stretched to all of
$\R^N$, the restriction map
\begin{align*}
  \Phi(\R^N) \to \Phi(0 \in \R^N)
\end{align*}
to the ``space'' of germs near 0 is an equivalence.  Now, a germ of a
graph around a point is easy to understand: Either it is the empty
germ, or it is the germ of a line through the point, or it is the germ
of $k\geq 3$ half-lines meeting at the point.  Any non-empty germ is
essentially determined by $k\geq 2$ points on $S^{N-1}$, so the space
of non-empty germs of graphs is essentially the space of finite
subsets of cardinality $\geq 2$ of $S^{N-1}$.  Let $\Sub(S^{N-1})$
denote the space of non-empty finite subsets of $S^{N-1}$.  The space
$\Sub(S^{N-1})$ is not quite right, for two reasons---it doesn't model
the empty germ, and it includes points that it shouldn't, namely the
space of 1-point subsets $S^{N-1} \subseteq \Sub(S^{N-1})$.  Both of
these problems can be fixed by collapsing the space of 1-point subsets
$S^{N-1} \subseteq \Sub(S^{N-1})$ to a point.  The above discussion
defines a map
\begin{align}\label{eq:13}
  \Phi(0 \in \R^N) \to \Sub(S^{N-1}) / S^{N-1},
\end{align}
which maps the empty germ to $[S^{N-1}]$ and maps the germ of $(G,0)$
to the set of tangent directions of $G$ at 0.  It seems reasonable
that this map should be a homotopy equivalence (it even seems close to
being a homeomorphism: If we had considered instead piecewise linear
graphs, it would be a bijection).  Curtis and To Nhu \cite{MR794488}
proves that $\Sub(S^{N-1})$ is contractible. (In fact they prove that
it is homeomorphic to $\R^\infty$.  For an easy, and more relevant,
proof of weak contractibility see \cite{MR1823966} or \cite[\S
3.4.1]{MR2058353}.)  Therefore the right hand side of the
map~(\ref{eq:13}) is homotopy equivalent to $S^N$ as we want.
Unfortunately, the natural map from $\Phi(\R^N)$ to the right hand
side of~\eqref{eq:13}, which assigns to $G \in \Phi(\R^N)$ the set of
directions of half-edges through $0 \in \R^N$, is not even continuous.

Let $\mathscr{D}$ be the category of finite sets and surjections.
Then
\begin{align*}
  \Sub(S^{N-1}) = \colim_{T \in \mathscr{D^\mathrm{op}}} \prod_T S^{N-1}.
\end{align*}
A step towards rectifying~(\ref{eq:13}) to a continuous map is to
replace the colimit by the homotopy colimit.  But the real reason for
discontinuity is that from the point of view of germs at a point, the
collapse of an edge leads to a \emph{sudden} splitting of one
half-edge into two.  To fix this, we will fatten up $\Phi(\R^N)$ in a
way that allows us to remove the suddenness of edge collapses, and
remotely similar to the proof of the equivalence~(\ref{eq:9}) in
section~\ref{subsec:poset-model-graph}.

\subsection{A pushout diagram}
\label{subsec:an-honest-proof}

The main result of this section is proposition~\ref{propcor:rephrased}
below.  Recall that a graph $G$ is a tree if it is contractible (in
particular non-empty).

\begin{definition}
  Let $\CC$ be the topological category whose objects are triples
  $(G,r,\phi)$, where $G \in \Phi(\R^N)$ and $r > 0$ satisfies that $G
  \pitchfork \partial B(0,r)$ and that $G \cap B(0,r)$ is a tree.
  $\phi$ is a labelling of the set of leaves, i.e.\ a bijection
  \begin{align*}
    \phi: \underline{k} = \{1, \dots, k\}
    \stackrel{\cong}{\longrightarrow} G \cap \partial B(0,r).
  \end{align*}
  Topologize $\Ob(\CC)$ as a subset
  \begin{align*}
    \Ob(\CC) \subseteq U_1 \times \coprod_{\substack{r>0\\k\geq 2}}
    \Map(\underline{k}, \partial B(0,r)).
  \end{align*}
  There is a unique morphism $(G,r,\phi) \to (G', r', \phi')$ if and
  only if $G = G'$ and $r\leq r'$, otherwise there is none.
\end{definition}

\begin{definition}\label{def:E-bullet}
  Let $E_\bullet$ be the simplicial space where an element of $E_k
  \subseteq N_k\CC \times S^N$ is a pair $(\chi,p)$, where $\chi =
  (G,r_0 < r_1 < \dots < r_k, \{\phi_i\}) \in N_k\CC$ and $p \in S^N =
  \R^N \cup \{\infty\}$ satisfies 
  \begin{align*}
    p \in \R^N \cup \{\infty\} - \big(G \cap B(0,r_k) - \Int
    B(0,r_0)\big).
  \end{align*}
  Include $N_\bullet\CC \subset E_\bullet$ as the subset with $p =
  \infty$.
\end{definition}
\begin{proposition}\label{propcor:rephrased}
  Let $B\CC \to |E_\bullet|$ be included as the subspace with $p =
  \infty$.  Then we have a weak equivalence $\Phi(\R^N) \simeq
  |E_\bullet|/B\CC$.
\end{proposition}

Proposition~\ref{propcor:rephrased} is proved in several steps.
First, in lemma~\ref{lemma:ho-pushout}, we write $\Phi(\R^N)$ as a
homotopy pushout of three open subsets $U_0$, $U_1$, and $U_{01}$.  In
lemma~\ref{lemma:ho-pushout-E} we give a similar description of
$|E_\bullet|/B\CC$ as a homotopy pushout.  Then we relate the homotopy
pushout diagrams by a zig-zag of weak equivalences maps according to
diagram~\eqref{eq:27}.
\begin{definition}\ 
  \begin{enumerate}[(i)]
  \item Let $U_0\subseteq \Phi(\R^N)$ be the subset consisting of
    graphs $G$ satisfying $0 \not \in G$.
  \item Let $U_1\subseteq \Phi(\R^N)$ be the subset consisting of
    graphs $G$ for which there exists an $r > 0$ such that $G
    \pitchfork \partial B(0,r)$ and that $G \cap B(0,r)$ is a tree.
  \item Let $U_{01} = U_0 \cap U_1$.
  \end{enumerate}
\end{definition}
\begin{lemma}\label{lemma:ho-pushout}
  The homotopy pushout (double mapping cylinder) of the diagram
  \begin{align}\label{eq:25}
    U_0 \leftarrow U_{01} \rightarrow U_1
  \end{align}
  is weakly equivalent to $\Phi(\R^N)$.
\end{lemma}
\begin{proof}
  $\Phi(\R^N)$ is the union of the two subsets $U_0$ and $U_1$, and it
  is easy to see that both of these are open.
\end{proof}
To give a similar pushout description of $|E_\bullet|/B\CC$ in
lemma~\ref{lemma:ho-pushout-E} below we need the following
definitions.
\begin{definition}
  Let $F_0$, $F_{01}$, and $F_1$ be the functors $\CC \to
  \mathrm{Spaces}$ given by
  \begin{align*}
    F_0(G,r,\phi) & = \R^N \cup\{\infty\} - G \cap B(0,r)\\
    F_{01}(G,r,\phi) & = \Int B(0,r) - G\\
    F_1(G,r,\phi) & = \Int B(0,r).
  \end{align*}
\end{definition}
$F_0$ is contravariant and $F_{01}$ and $F_1$ are covariant.  All
three spaces $N_k(F_0 \wr \CC)$, $N_k(\CC \wr F_{01})$ and $N_k(\CC
\wr F_1)$ are open subsets of $N_k\CC \times S^N$, where $S^N = \R^N
\cup \{\infty\}$.
\begin{lemma}
  \label{lemma:ho-pushout-E}
  $|E_\bullet|$ is weakly equivalent to the homotopy pushout of the
  diagram
  \begin{align}\label{eq:35}
    B(F_0 \wr\CC) \leftarrow B(\CC\wr F_{01}) \rightarrow B(\CC
    \wr F_1),
  \end{align}
  and $|E_\bullet|/B\CC$ is weakly equivalent to the homotopy pushout
  of the diagram
  \begin{align}\label{eq:43}
    B(F_0 \wr\CC)/B\CC \leftarrow B(\CC\wr F_{01}) \rightarrow B(\CC
    \wr F_1).
  \end{align}
\end{lemma}
\begin{proof}
  As subsets of $N_k\CC \times S^N$ we have
  \begin{align*}
    N_k(F_0 \wr \CC) \cap N_k(\CC \wr F_1) & = N_k(\CC \wr F_{01})\\
    N_k(F_0 \wr \CC) \cup N_k(\CC \wr F_1) & = E_k.
  \end{align*}
  Then $E_k$ is weakly equivalent to the homotopy pushout of the
  following diagram
  \begin{align}\label{eq:45}
    N_k(F_0 \wr \CC) \leftarrow N_k(\CC \wr F_{01}) \rightarrow
    N_k(\CC \wr F_1).
  \end{align}
  But the homotopy pushout of diagram~\eqref{eq:35} is homeomorphic to
  the geometric realization of the simplicial space whose
  $k$-simplices is the homotopy pushout of~\eqref{eq:45}.  The second
  part is similar.
\end{proof}

We will now relate the pushout diagram~\eqref{eq:25} to the pushout
diagram~\eqref{eq:43} by a zig-zag of maps, according to the following
diagram.
\begin{align}\label{eq:27}
\begin{aligned}
  \xymatrix{
    U_0 & U_{01} \ar[l]\ar[r] & U_1 \\
    U_0 \ar@{=}[u]\ar[d] & {B\CC_{01}} \ar[l]\ar[u]\ar[r]\ar[d] & {B\CC}
    \ar[u] \ar[d]\\
    {*} & {B(\CC\wreath F_{01})} \ar[l]\ar[r] & {B(\CC \wreath F_1)}\\
    {B(F_0\wreath\CC)/B\CC} \ar[u] & {B(\CC\wreath
      F_{01})}\ar[l]\ar@{=}[u] \ar[r] & {B(\CC \wreath F_1).} \ar@{=}[u] }
\end{aligned}
\end{align}
The spaces and maps in the diagram will be defined below, and we will
prove that all vertical maps are weak equivalences.  We first consider
the second row of the diagram.
\begin{definition}
  Let $\CC_{01}$ be the subcategory of $\CC$ consisting of $(G,r,\phi)$
  with $G \in U_{01}$.
\end{definition}
\begin{proposition}\label{proposition:U1}
  The forgetful maps
  \begin{align*}
    B\CC \to U_1,\quad B\CC_{01} \to U_{01}
  \end{align*}
  are both weak equivalences.
\end{proposition}
\begin{proof}
  This is completely similar to theorem \ref{lemthm:equivalences} and
  lemma~\ref{lemma:empty-intersections}: $N_k\CC \to U_1$ is etale for
  all $k$, and for $G \in U_1$, the inverse image in $\CC$ is
  equivalent as a category to a totally ordered non-empty set.
  Similarly for $B\CC_{01} \to U_{01}$.
\end{proof}

Maps from the second to the third row in diagram~\eqref{eq:27} are
induced by the natural diagram of functors
\begin{align}\label{eq:14}
\begin{aligned}
  \xymatrix{
    {\CC_{01}} \ar[d]\ar[r] & {\CC}\ar[d]\\
    {\CC \wreath F_{01}} \ar[r]
    & {\CC \wreath F_1}.
  }
\end{aligned}
\end{align}
The horizontal functors in~\eqref{eq:14} are the natural inclusions,
and the vertical functors are both given by $(G,r,\phi) \mapsto
(G,r,\phi,0)$.

\begin{lemma}\ 
  \begin{enumerate}[(i)]
  \item $U_0$ is contractible.
  \item $B\CC \to B(\CC \wreath F_{1})$ is a weak equivalence.
  \item $B\CC_{01} \to B(\CC\wreath F_{01})$ is a weak equivalence.
  \end{enumerate}
\end{lemma}
\begin{proof}
  (i) follows by pushing radially away from $p=0$, as in
  lemma~\ref{lemma:Phi-is-connected}.

  (ii) is also easy.  Moving the point $p \in \Int B(0,r)$ to 0 along
  a straight line defines a deformation retraction of $B(\CC \wreath
  F_1)$ onto the image of $B\CC$.

  For (iii), notice that for each $k$ we have the following pullback
  diagram of spaces
  \begin{align*}
    \xymatrix{
      {N_k \CC_{01}} \ar[r]\ar[d] & {N_k (\CC\wreath F_{01})} \ar[d] \\
      {\displaystyle\coprod_{r > 0} \{0\}} \ar[r] & {
        \displaystyle\coprod_{r > 0} \Int B(0,r)}.  }
  \end{align*}
  It is easy to see that the right hand vertical map is a fibration
  (in fact a trivial fiber bundle), so the diagram is also homotopy
  pullback.  The bottom horizontal map is obviously a homotopy
  equivalence, so it follows that $N_k \CC_{01} \to N_k(\CC\wreath
  F_{01})$ is an equivalence for all $k$.  This proves (iii).
\end{proof}

The map from the third to the fourth row of diagram~\eqref{eq:27} is
covered by the following lemma.
\begin{lemma}
  The inclusion $\{\infty\} \to F_0(G,r,\phi)$ is a homotopy
  equivalence, and $B(F_0 \wr\CC)/B\CC$ is weakly contractible.
\end{lemma}
\begin{proof}
  $N_k(F_0 \wr \CC)$ is an open subset of $N_k\CC \times S^N$ such
  that all fibers of the projection
  \begin{align*}
    N_k(F_0 \wr \CC) \to N_k\CC
  \end{align*}
  are contractible.  It follows from \cite[proposition
  (A.1)]{MR516216} that the projection is a Serre fibration and hence
  a weak equivalence.  Therefore the section $N_k\CC \to N_k(F_0\wr
  \CC)$ obtained by setting $p = \infty$ is also a weak equivalence.
  It is easy to see that this section is a cofibration, so the
  quotient
  \begin{align*}
    N_k(F_0 \wr \CC)/N_k\CC
  \end{align*}
  is weakly contractible.
\end{proof}
This finishes the proof of proposition~\ref{propcor:rephrased}.

\begin{remark}\label{remark:generalized-spherical-fib}
  From the third line in diagram~\eqref{eq:27} it follows that
  $\Phi(\R^N)$ is weakly equivalent to the mapping cone of the map
  $B(\CC\wreath F_{01}) \to B\CC$.  One can think of this map as a
  ``generalized spherical fibration'', and hence of the mapping cone
  as a ``generalized Thom space'', in the following sense.  The fiber
  of the map
  \begin{align*}
    N_k(\CC\wreath F_{01}) \to N_k \CC
  \end{align*}
  over a point $(G, r_0 < r_1 < \dots < r_k, \{\phi_i\})$ is the space
  \begin{align*}
    \Int B(0,r_0) - G \simeq \bigvee^{k_0-1} S^{N-2},
  \end{align*}
  where $k_0$ is the cardinality of the set $G \cap \partial
  B(0,r_0)$.  Thus, the fibers of $B(\CC\wreath F_{01}) \to B\CC$ are
  not spheres, as they would be were the map an honest spherical
  fibration, but wedges of spheres, where the number of spheres in the
  fiber varies over the base.
\end{remark}

\subsection{A homotopy colimit decomposition}
\label{subsec:alexander-duality}

Let $\DD_{\geq 2}$ be the category whose objects are finite sets of
cardinality at least 2, and whose morphisms are the surjective maps of
sets.  In this section we will first rewrite $|E_\bullet|/B\CC$ stably
as the pointed homotopy colimit of a functor $H: \DD_{\geq
  2}^\mathrm{op} \to \mathrm{Spaces}$.  This is done in
proposition~\ref{prop:E-is-hocolim} below.  Then we prove that this
pointed homotopy colimit is weakly equivalent to $S^N$ in
proposition~\ref{prop:spaces-contractible}.  Together these results
prove theorem~\ref{thm:homotopy-type-Phi}.

There is a functor $T: \CC \to \DD_{\geq 2}^\mathrm{op}$ defined in
the following way.  Let $(G, r, \phi) \to (G, r', \phi')$ be a
morphism in $\CC$.  We have a diagram of inclusions
\begin{align*}
  \xymatrix{
    {G \cap \partial B(0,r)}\ar[r]^-{i_r} & {G \cap (B(0,r') -
      \Int B(0,r))} & {G \cap \partial B(0,r')}\ar[l]_-{i_{r'}}
  }
\end{align*}
in which the inclusion $i_r$ is a homotopy equivalence and $i_{r'}$
induces a surjection in $\pi_0$.
\begin{definition}\label{defn:eq-26}
  Let $f: (G, r, \phi) \to (G, r', \phi')$ be a morphism, and let
  $i_r$ and $i_{r'}$ be as above.  Then let $T(f)$ be the composition
  \begin{align*}
    \phi^{-1}\circ(\pi_0 i_r)^{-1} \circ (\pi_0 i_{r'}) \circ \phi':
    \underline{k}' \to \underline{k}.
  \end{align*}
  This defines a functor $T: \CC \to \DD_{\geq 2}^\mathrm{op}$.
\end{definition}
\begin{definition}
  For $\underline k \in \DD_{\geq 2}^\mathrm{op}$, let $\Delta
  \subseteq (S^{N-1})^{\underline k} = \Map(\underline k, S^{N-1})$ be
  the diagonal.  Let the functor $H: \DD_{\geq 2}^\mathrm{op} \to
  \mathrm{Spaces}$ be the quotient
  \begin{align*}
    H(\underline{k}) = \Map(\underline{k},S^{N-1})/\Delta.
  \end{align*}
\end{definition}
The following proposition will be proved below in several steps.  We
will say that a map is ``highly connected'' if it is $c(N)$-connected
for a function $c: \N \to \N$ such that $c(N) \to \infty$ as $N \to
\infty$.  Similarly we will say that a map is ``$N+$highly connected''
if it is $(N + c(N))$-connected.
\begin{proposition}
  \label{prop:E-is-hocolim}
  There is an $N+$highly connected map
  \begin{align*}
    |E_\bullet|/B\CC \to B(\DD_{\geq 2}^\mathrm{op} \wr H)/B\DD_{\geq
      2}^\mathrm{op}.
  \end{align*}
\end{proposition}
Recall that the space $B(\DD_{\geq 2}^\mathrm{op} \wr H)$ is the
homotopy colimit of $H$.  Each $H(\underline k)$ has the basepoint
$[\Delta]$ which defines an inclusion $B\DD_{\geq 2}^\mathrm{op}
\subset B(\DD_{\geq 2}^\mathrm{op} \wr H)$.  The quotient space is the
\emph{pointed} homotopy colimit of the functor $H$.

Let $K\subseteq \R^N$ be a compact subset with contractible path
components.  The \emph{duality} map is the map
\begin{align*}
  A: (\R^N - K) \to \Map(K, S^{N-1})
\end{align*}
given by
\begin{align*}
  A(p)(x) = \frac{p-x}{|p-x|}
\end{align*}
The map $A$ is $(2N-3)$-connected.  Indeed, it is homotopy equivalent
to the inclusion
\begin{align*}
  \bigvee^{\pi_0K} S^{N-1} \to \prod_{\pi_0K} S^{N-1}.
\end{align*}
Let $\Delta \subseteq \Map(K,S^{N-1})$ denote the constant maps.  $A$
induces a well defined, continuous map
\begin{align}\label{eq:44}
  \R^N \cup\{\infty\} - K \xrightarrow{ A } \Map(K,S^{N-1}) / \Delta
\end{align}
by mapping $\infty\mapsto [\Delta]$.  This map is also
$(2N-3)$-connected.

As $K$, we can take the space $G \cap B(0,r_k) - \Int B(0,r_0)$ in the
definition of $E_\bullet$.  This leads to the following definition.
\begin{definition}\label{defn:tildeE}
  Let $\tilde E_\bullet$ be the simplicial space where an element of
  $\tilde E_k$ is a pair $(\chi,f)$, where $\chi = (G,r_0 < r_1 <
  \dots < r_k, \{\phi_i\}) \in N_k\CC$ and $f$ is an element
  \begin{align*}
    f \in \Map(K,S^{N-1})/\Delta,
  \end{align*}
  where $K = G \cap B(0,r_k) - \Int B(0,r_0)$ and $\Delta$ denotes the
  subset of constant maps.
\end{definition}
The subset $K$ in the above definition will be a \emph{forest}, i.e.\
a disjoint union of (at least two) contractible graphs.  We should
explain the topology on the space $\tilde E_k$.  The main observation
is that if $\chi,\chi' \in N_k\CC$ and $K,K'$ are the corresponding
forests, then there will be a canonical map $\phi: K' \to K$ whenever
$\chi'$ is sufficiently close to $\chi$.  (By the definition of the
topology on $\Phi(\R^N)$, any $G'$ near $G$ will admit a map $\tilde
\phi: G' \dashrightarrow G$ whose domain contains $K'$ and whose image
contains $K$.  After reparametrizing edges it will restrict to a map
from $K'$ onto $K$.)  We topologize $\tilde E_k$ by declaring
$(\chi',f')$ close to $(\chi,f)$ if $\chi'$ is close to $\chi$ and
$f'$ is close to $f \circ \phi$.

For the following lemma, recall the notion of \emph{fiber homotopy}
from \cite{MR0155330}, and some related notions.  If $f:E \to B$ and
$f':E' \to B$ are two maps, then a \emph{fiber homotopy} is a homotopy
$F: E \times [0,1] \to E'$ over $B$.  A map $g: E \to E'$ over $B$ is
a \emph{fiber homotopy equivalence} if it admits a map $h: E' \to E$
which is left and right inverse to $g$ up to fiber homotopy.  A map $E
\to B$ is \emph{fiber homotopy trivial} if it is fiber homotopy
equivalent to a projection $B \times F \to B$.  A map $f: E \to B$ is
\emph{locally fiber homotopy trivial} if $B$ admits a covering by open
sets $U$ such that the restriction $f^{-1}(U) \to U$ is fiber homotopy
trivial.  It is shown in \cite[theorem 6.4]{MR0155330} that local
fiber homotopy triviality is sufficient for the ``long exact sequence
for a fibration'': if $f:E \to B$ is locally fiber homotopy trivial,
then the homotopy groups of a fibers $F_b = f^{-1}(b)$ fit into a long
exact sequence with $\pi_*(E)$ and $\pi_*(B)$.

It follows from the definition that the projection $\tilde E_k \to N_k
\CC$ is locally fiber homotopy trivial.  Indeed, let $U \subseteq
N_k\CC$ be a neighborhood of $\chi$ small enough that any $\chi' \in
U$ admits a canonical map $\phi: K' \to K$ (cf.\ the discussion
following definition~\ref{defn:tildeE}).  We get a map
\begin{align*}
  U \times \Map(K,S^{N-1}) \to \tilde E_k,
\end{align*}
given by $(\chi',f)\mapsto (\chi',f \circ \phi)$, which restricts to a
fiber homotopy equivalence over $U$.

\begin{lemma}\label{lemcor:map-induced-by-A}
  The map $A$ above induces $N+$highly connected maps $|E_\bullet| \to
  |\tilde E_\bullet|$ and $|E_\bullet|/B\CC \to |\tilde
  E_\bullet|/B\CC$.
\end{lemma}
\begin{proof}
  Both maps $E_k \to N_k\CC$ and $\tilde E_k \to N_k\CC$ induce long
  exact sequences in homotopy groups.  For $\tilde E_k$, this was
  explained above, and for $E_k$ it can be proved in the following
  way.  We shall prove later (lemma~\ref{lemma:C-approx-D}) that
  $N_k\CC$ has trivial $\pi_1$ and $\pi_2$.  A similar argument shows
  that $E_k$ is simply connected.  Hence the homotopy fiber of the
  projection $E_k \to N_k\CC$ is simply connected.  From
  \cite[proposition 5]{MR0402733} it follows that the inclusion of the
  fiber of the projection $E_k \to N_k\CC$ into the homotopy fiber is
  a \emph{homology equivalence}.  Since both the homotopy fiber and
  the fiber are simply connected, it is actually a homotopy
  equivalence, so $E_k \to N_k\CC$ is a quasifibration.

  The induced map on fibers is the map~\eqref{eq:44}, so the first
  part of the lemma follows from the 5-lemma.  The second map uses
  that the inclusions of $B\CC$ into $|E_\bullet|$ and $|\tilde
  E_\bullet|$ are cofibrations.
\end{proof}

An element of $N_l(\CC \wreath (H\circ T))$, where $T$ is the functor
from definition~\ref{defn:eq-26}, is given by an element $(G, r_0 <
r_1 < \dots < r_l,\{\phi_i\}) \in N_l(\CC)$ together with an element
$g\in\Map(\underline{k_0},S^{N-1})/\Delta$.  Here,
\begin{align*}
  \phi_i : \underline{k_i} \to G \cap \partial B(0,r_i)
\end{align*}
are the labellings.  Again, let $K = G \cap B(0,r_l) - \Int B(0,r_0)$.
The labelling $\phi_0$ in the first vertex induces an injective map
\begin{align*}
  \phi_0: \underline{k_0} \to K
\end{align*}
which is a homotopy equivalence.  It has a unique left inverse which
we denote $\phi_0^{-1}$.  Up to homotopy $\phi_0^{-1}$ is also right
inverse to $\phi_0$.

Composition with $\phi_0^{-1}$ induces a homotopy equivalence
\begin{align*}
  \Map(k_0, S^{N-1})/\Delta \xrightarrow{ \circ \phi_0^{-1} } \Map(K,
  S^{N-1})/\Delta
\end{align*}
and in turn a simplicial map
\begin{align}
  \label{eq:28}
  N_\bullet(\CC \wreath (H\circ T)) \to \tilde E_\bullet
\end{align}
which is a degreewise homotopy equivalence.  Similarly to
lemma~\ref{lemcor:map-induced-by-A}, this proves the following lemma.
\begin{lemma}\label{lemcor:H-circ-T}
  The maps
  \begin{align*}
    B(\CC \wreath (H\circ T)) & \to |\tilde E_\bullet|,\\
    B(\CC \wreath (H\circ T))/B\CC & \to |\tilde E_\bullet|/B\CC
  \end{align*}
  induced by~\eqref{eq:28} are weak homotopy equivalences.\qed
\end{lemma}

Combining proposition~\ref{propcor:rephrased} and
lemmas~\ref{lemcor:map-induced-by-A} and \ref{lemcor:H-circ-T}, we
get the following.
\begin{corollary}\label{cor:19}
  There is an $N+$highly connected map $\Phi(\R^N) \to B(\CC \wreath
  (H \circ T))/B\CC$.\qed
\end{corollary}

Corollary~\ref{cor:19} states that stably (i.e.\ for $N \to \infty$),
we can regard $\Phi(\R^N)$ as the pointed homotopy colimit of the
functor $(H\circ T)$ over the topological category $\CC$.  We would
like to replace that with the pointed homotopy colimit of the functor
$H$ over the category $\DD_{\geq 2}$, whose objects are finite sets
$\underline{k}$ of cardinality at least 2 and whose morphisms are
surjections.
\begin{lemma}\label{lemma:C-approx-D}
  The functor $T: \CC \to \DD_{\geq 2}^\mathrm{op}$ induces a highly
  connected map $N_lT: N_l \CC \to N_l \DD_{\geq 2}^\mathrm{op}$ for
  all $l$.
\end{lemma}
\begin{proof}
  The codomain $N_l \DD_{\geq 2}$ is a discrete set.  Let
  $(\underline{k_0} \to \underline{k_1} \to \dots \to \underline{k_l})
  \in N_l\DD^\mathrm{op}_{\geq 2}$.  A point in the inverse image is
  given by embeddings of the finite sets $\underline{k_i}$ into
  $(N-1)$-spheres, and trees with these sets as the set of leaves.
  Embeddings of finite sets into an $(N-1)$-sphere form an
  $(N-3)$-connected space.  Trees with a fixed set of leaves form an
  $(N-4)$-connected space by theorem~\ref{thm:conclusion-htpy-types}
  (applied with $M = B(0,a_j) - \Int B(0,a_{j-1})$, and using that
  $A_0^s$ is the trivial group).
\end{proof}

The approximation in lemma~\ref{lemma:C-approx-D} may seem to be not
good enough.  $\Omega^\infty\mathbf{\Phi}$ is the direct limit of the
spaces $\Omega^N\Phi(\R^N)$, so we should deal with spaces up to
$N+$highly connected maps instead of just up to highly connected maps.
Surprisingly, the extra $N$ comes for free.  (Analogously, if $f: X
\to Y$ is $c$-connected and $\xi$ is an $N$-dimensional vector bundle
over $Y$, then the map of Thom spaces $X^{f^*\xi} \to Y^\xi$ is
$(c+N)$-connected.)  Proposition~\ref{prop:asdf} finishes the proof of
proposition~\ref{prop:E-is-hocolim}.
\begin{proposition}\label{prop:asdf}
  The map
  \begin{align*}
    B(\CC\wreath (H\circ T))/B\CC \to B(\DD_{\geq
      2}^\mathrm{op}\wreath H) /B\DD_{\geq 2}^\mathrm{op}
  \end{align*}
  is $N+$highly connected.
\end{proposition}
\begin{proof}
  For all $k$ we have the following pullback diagram.
  \begin{align*}
    \xymatrix{ {N_k(\CC\wreath (H\circ T))}\ar[r]\ar[d] &
      {N_k(\DD_{\geq 2}^\mathrm{op} \wreath
        H)} \ar[d]\\
      {N_k \CC} \ar[r] & {N_k \DD_{\geq 2}^\mathrm{op}}.  }
  \end{align*}
  The right hand vertical map is a fibration, so the diagram is also
  homotopy cartesian.  Both vertical maps are split, using the
  canonical basepoint $\infty \in H$.  It follows that the diagram
  \begin{align*}
    \xymatrix{
      {N_k \CC} \ar[r]\ar[d] & {N_k \DD_{\geq 2}^\mathrm{op}}\ar[d]\\
      {N_k(\CC\wreath (H\circ T))}\ar[r] & {N_k(\DD_{\geq
          2}^\mathrm{op} \wreath H)} }
  \end{align*}
  is also homotopy cartesian (horizontal homotopy fibers are homotopy
  equivalent).

  The vertical and horizontal maps are all $(N-3)$-connected.  It
  follows by the Blakers-Massey theorem that the diagram is $(N-3) +
  (N-3) - 1 = (2N-7)$-cocartesian.  This means precisely that the
  induced map of vertical cofibers is $(2N-7)$-connected and the claim
  follows.
\end{proof}
Thus, we have an $N+$highly connected map from $\Phi(\R^N)$ to the
pointed homotopy colimit of the functor $H: \DD_{\geq 2}^\mathrm{op}
\to \mathrm{Spaces}$.  We proceed to determine the homotopy type of
this pointed homotopy colimit.  Recall that $H(\underline{k}) =
\Map(\underline{k},S^{N-1})/\Delta$.  The pointed homotopy colimit is
homeomorphic to the quotient
\begin{align}\label{eq:15}
  B(\DD_{\geq 2}^\mathrm{op}\wreath \Map(-,S^{N-1}))/B(\DD_{\geq
    2}^\mathrm{op}\wreath \Delta)
\end{align}
where $\Delta$ denotes the constant functor $S^{N-1}$.
\begin{proposition}\label{prop:spaces-contractible}
  The spaces $B(\DD_{\geq 2}^\mathrm{op}\wreath\Map(-,S^{N-1}))$ and
  $B\DD_{\geq 2}$ are both contractible.
\end{proposition}
\begin{proof}[Proof of theorem~\ref{thm:homotopy-type-Phi}]
  Proposition \ref{prop:spaces-contractible} implies that $B(\DD_{\geq
    2}^\mathrm{op}\wreath \Delta) \cong B\DD_{\geq 2}^\mathrm{op}
  \times S^{N-1} \simeq S^{N-1}$, so the quotient in~(\ref{eq:15})
  becomes $S^N$ and we get an $N+$highly connected map
  \begin{align}
    \label{eq:57}
    \Phi(\R^N) \to S^N
  \end{align}
  The map~(\ref{eq:57}) is a zig-zag of $N+$highly connected maps, all
  of which induce spectrum maps as $N$ varies.  It follows that there
  is a weak equivalence of spectra $\mathbf{\Phi} \simeq S^0$ as
  claimed.
\end{proof}

\begin{remark}
  For an object $\underline{k} \in \DD_{\geq 2}$, let $\Delta \to
  (S^{-1})^{\underline{k}}$ be the inclusion of the diagonal into
  the $k$-fold power of the spectrum $S^{-1}$.  Let
  $(S^{-1})^{\underline{k}}/\Delta$ be the cofiber.  Then we have
  proved two homotopy equivalences
  \begin{align*}
    \mathbf{\Phi} \simeq \hocolim_{\underline{k} \in \DD_{\geq 2}}
    \bigg((S^{-1})^{\underline{k}}/\Delta\bigg) \simeq S^0.
  \end{align*}
\end{remark}

\begin{proof}[Proof of proposition~\ref{prop:spaces-contractible}]
  We have a functor $\DD_{\geq 2} \to \DD_{\geq 2}$ given by $T
  \mapsto \underline{2} \times T$, and the projections define natural
  transformations
  \begin{align*}
    \xymatrix{
      T & {\underline{2} \times T} \ar[l]\ar[r] & {\underline{2}}.
    }
  \end{align*}
  This contracts $B\DD_{\geq 2}$ to the point $\underline{2} \in B \DD_{\geq 2}$.

  For the space
  \begin{align*}
    B(\DD_{\geq 2}^\mathrm{op}\wreath\Map(-,S^{N-1})) =
    \hocolim_{T \in \DD_{\geq 2}^\mathrm{op}} \Map(T,S^{N-1})
  \end{align*}
  we use a trick strongly inspired by the works of \cite{MR1823966}
  and \cite[\S 3.4.1]{MR2058353}, which prove that the colimit (not
  homotopy colimit) is contractible.

  Choose a (symmetric monoidal) disjoint union functor $\amalg:
  \DD_{\geq 2}\times \DD_{\geq 2} \to \DD_{\geq 2}$.  For brevity,
  denote the functor $\Map(-,S^{N-1})$ by $J$.  The disjoint union
  functor induces a functor
  \begin{align*}
    (\DD_{\geq 2}^\mathrm{op} \wreath J) \times (\DD_{\geq
      2}^\mathrm{op} \wreath J) \to
    (\DD_{\geq 2}^\mathrm{op}\wreath J)
  \end{align*}
  which is associative and commutative up to natural transformation.
  It follows that the classifying space is a homotopy associative and
  homotopy commutative $H$-space.

  In this $H$-space structure, multiplication by 2 is homotopic to the
  identity.  This follows from the natural transformation $T \amalg T
  \to T$.  The claim then follows from lemma~\ref{lemma:mult-by-2}
  below.
\end{proof}
\begin{lemma}\label{lemma:mult-by-2}
  A connected, homotopy associative, homotopy commutative $H$-space
  $X$ is weakly contractible if multiplication by 2 (i.e.\ the map $x
  \mapsto x \cdot x$) is homotopic to the identity.
\end{lemma}
This lemma is completely trivial when $X$ has a homotopy unit.  In
that case, it is well know that the map induced by the $H$-space
structure
\begin{align*}
  \pi_*X \times \pi_*X \to \pi_*X
\end{align*}
agrees with the usual group multiplication on homotopy groups.  Hence
all $x \in \pi_*X$ satisfies $x + x = x$.  The proof in the general
case is a variation of this argument.
\begin{proof}
  Let $\mu: X \times X \to X$ be the $H$-space structure.  Choose a
  basepoint $x_0 \in X$ and write $\pi_n(X) = \pi_n(X,x_0)$.  We can
  assume that $\mu$ is a pointed map.  The two projections $p,q: X
  \times X \to X$ induce an isomorphism
  \begin{align*}
    (p_*,q_*): \pi_n(X \times X) \to \pi_nX\times \pi_nX
  \end{align*}
  and we let
  \begin{align*}
    \bullet = \mu_* \circ (p_*,q_*)^{-1}: \pi_nX \times \pi_nX \to
    \pi_nX.
  \end{align*}
  This is now an associative, commutative product on $\pi_nX$
  satisfying $x \bullet x = x$ for all $x \in \pi_n(X)$.  Let $+$
  denote the usual group structure on $\pi_nX$ and write $0$ for the
  identity element with respect to $+$ (we will write it additively
  although we don't yet know that it is commutative for $n = 1$).

  Let $\Delta: X \to X \times X$ be the diagonal and $i,j: X \to X
  \times X$ the inclusions $i(x) = (x,x_0)$, $j(x) = (x_0,x)$.  Then
  we have
  \begin{align*}
    (p_*,q_*) \circ (i_* + j_*)(x) &= ((p\circ i)_*, (q\circ i)_*)(x) +
    ((p\circ j)_*, (q\circ j)_*)(x) \\ &= (x,0) + (0,x) = (x,x) =
    (p_*,q_*)\circ \Delta_*(x).
  \end{align*}
  It follows that $\Delta_* = i_* + j_*$ because $(p_*,q_*)$ is an
  isomorphism.  Now $\mu\circ \Delta \simeq \id$ and $\mu \circ i
  \simeq \mu \circ j$, so
  \begin{align*}
    x = \mu_*\Delta_*(x) = \mu_*i_*x + \mu_*j_*x = 2 x \bullet x_0.
  \end{align*}
  Substituting $x\bullet x_0$ for $x$ then gives
  \begin{align*}
    x \bullet x_0 = 2 x \bullet x_0 \bullet x_0 = 2 x \bullet x_0 = x
  \end{align*}
  which in turn gives that $x = 2x$ for any $x \in \pi_nX$.  It
  follows that $X$ is weakly contractible.
\end{proof}

\section{Remarks on manifolds}
\label{sec:some-remarks-manif}

Most of the results of this paper works equally well for the sheaf
$\Psi_d$, where $\Psi_d(U)$ is the space of all closed sets $M\subseteq
U$ which are smooth $d$-dimensional submanifolds without boundary.  A
neighborhood basis at $M$ is formed by the sets
\begin{align*}
  \mathcal{V}_{K,W} = \{N \in \Psi_d(U) | \text{$N \cap K = j(M) \cap
    K$ for some $j \in W$ } \},
\end{align*}
where $K\subseteq U$ is a compact set and $W\subseteq \Emb(M,U)$ is a
neighborhood of the inclusion in the Whitney $C^\infty$ topology.

The analogues of propositions~\ref{prop:3A1} and \ref{prop:3A2} hold
with almost identical proofs and give the following weak equivalence.
\begin{align}
  \label{eq:1}
  B \mathcal{C}_d^N \simeq \Omega^{N-1} \Psi_d(\R^N)
\end{align}
Here $\mathcal{C}_d^N$ is the cobordism category whose objects are
closed $(d-1)$-manifolds $M \subseteq \{a\} \times \R^{N-1}$ and whose
morphisms are compact $d$-manifolds $W \subseteq [a_0, a_1] \times
\R^{N-1}$, cf.\ \cite[section 2]{math.AT/0605249}.  We note that for
the proof of lemma~\ref{lem:Bc-Bi-Bu}\ref{item:2} we will no longer
necessarily have $D_{N,k}$ connected; however it will be a grouplike
topological monoid, which suffices for the proof.

Let $\Gr_d(\R^N)$ be the Grassmannian of $d$-planes in $\R^N$, and
$U^\perp_{d,N}$ the canonical $(N-d)$-dimensional vector bundle over
it.  A point in $U^\perp_{d,N}$ is given by a pair $(V,v)\in
\Gr_d(\R^N)\times \R^N$ with $v \perp V$.  Let
\begin{align*}
  q: U_{d,N}^\perp \to \Psi_d(\R^N)
\end{align*}
be the map given by $q(V,v) = V-v \in \Psi_d(\R^N)$.  This gives a
homeomorphism onto the subspace of manifolds $M^d \subseteq \R^N$
which are affine subspaces.  $q$ extends continuously to the one-point
compactification of $U_{d,N}^\perp$ by letting $q(\infty) =
\emptyset$.  This one-point compactification is the \emph{Thom space}
$\thh(U_{d,N}^\perp)$, and we get a map
\begin{align}\label{eq:19}
  q: \thh(U_{d,N}^\perp) \to \Psi_d(\R^N).
\end{align}
We will show that~\eqref{eq:19} is a weak equivalence.  Define two
open subsets $U_0 \subseteq \Psi_d(\R^N)$ and $U_1 \subseteq
\Psi_d(\R^N)$ in the following way.  $U_0$ is the space of
$d$-manifolds $M$ such that $0 \not \in M$, and $U_1$ is the space of
manifolds such that the function $p \mapsto |p|^2$ has a unique,
non-degenerate minimum on $M$.  Let $U_{01} = U_0 \cap U_1$.  These
are open subsets, and $\Psi_d(\R^N)$ is the pushout of $(U_0
\leftarrow U_{01} \rightarrow U_1)$.
\begin{lemma}\label{lemma:6.1}
  Each restriction of $q$
  \begin{align*}
    q^{-1}(U_0) &\to U_0 \\
    q^{-1}(U_{01}) &\to U_{01} \\
    q^{-1}(U_1) &\to U_1     
  \end{align*}
  is a weak homotopy equivalence.  Consequently~(\ref{eq:19}) is a
  weak equivalence.
\end{lemma}
\begin{proof}
  $U_0$ and $q^{-1}(U_0)$ are both contractible: $q^{-1}(U_0)$
  contracts to the point $\infty$, and the path constructed in the
  proof of lemma~\ref{lemma:Phi-is-connected} gives a contraction of
  $U_0$, pushing everything to infinity, radially away from 0.

  For $U_1$ a deformation retraction is defined as follows.  Let $M
  \in U_1$ have $p$ as unique minimum of $p \mapsto |p|^2$.  Let
  $\phi_t(x) = p + (1-t)(x-p)$.  This defines a diffeomorphism $\R^N
  \to \R^N$ for $t < 1$.  A path $\gamma$ in $U_1$ from $M$ to a point
  in the image of $q$ is defined by $\gamma(t) = \phi_t^{-1}(M)$ for
  $t < 1$ and $\gamma(1) = p + T_pM$.  This proves that $q^{-1}(U_1)
  \to U_1$ is a deformation retraction.  This deformation restricts to
  a deformation retraction of $q^{-1}(U_{01}) \to U_{01}$.
\end{proof}
We have proved the following result.
\begin{proposition}\label{prop:6.2}
  $q: \thh(U_{d,N}^\perp) \to \Psi_d(\R^N)$ is a weak equivalence.
  \qed
\end{proposition}

Thus we have proved the following theorem.  In the limit $N \to
\infty$ we recover the main theorem of \cite{math.AT/0605249}, but
theorem~\ref{thm:GMTW-strengthening} holds also for finite $N$.
\begin{theorem}\label{thm:GMTW-strengthening}
  There is a weak homotopy equivalence
  \begin{align}\tag*{\qed}
    B\mathcal{C}_d^N \simeq \Omega^{N-1}\thh(U_{d,N}^\perp).
  \end{align}
\end{theorem}

\nocite{*}

\bibliographystyle{alpha}
\bibliography{graphs}

\end{document}